\documentclass{article}
\usepackage[latin1]{inputenc}
\usepackage{amsmath, amssymb, amsthm, latexsym, mathtools, mathrsfs}
\usepackage{mathrsfs,MnSymbol,graphicx}
\usepackage{color}
\usepackage[shortlabels]{enumitem}

\usepackage{tikz}
\usepackage{graphicx}

\newtheorem{theorem}{Theorem}[section]
\newtheorem{lemma}[theorem]{Lemma}
\newtheorem{proposition}[theorem]{Proposition}
\newtheorem{remark}[theorem]{Remark}

\newtheorem{defin}[theorem]{Definition}
\newtheorem{obs}[theorem]{Observation}

\def\vect{\overrightarrow}
\def\rr{{\mathbb R}}

\def\am{^{-1}}

\def\su{\subset}
\def\se{\setminus}

\def\Ga{\Gamma}

\def\stb{,\ldots ,}
\def\stp{+ \ldots +}
\def\noi{\noindent}
\def\ov{\overline}
\def\Om{\Omega}
\def\proof{\noi {\bf Proof.} }

\def\noi{\noindent}
\def\msk{\medskip}
\def\bsk{\bigskip}

\def\ep{\varepsilon}
\def\rg{\mathring}
\begin{document}
\title{Decomposition of balls in $\rr^d$}

\author{Gergely Kiss and G\'abor Somlai \thanks{Research supported by the Hungarian National Foundation for Scientific Research, Grant
No. T49786.}~ } \footnotetext[1]{{\bf Keywords:}
$m$-divisibility of the balls in Euclidean space}
\footnotetext[2]{{\bf MSR 2010 classification:} 52A15, 52A20, 51F20}
\footnotetext[3]{The research was partially
supported by the Hungarian National Foundation for Scientific Research,
Grant No. K104178}
\date{}
\maketitle

\begin{abstract}
We investigate the decomposition problem of balls into finitely many congruent pieces in dimension $d=2k$. In addition,
we prove that the $d$ dimensional unit ball $B_d$ can be divided into finitely many
congruent pieces if $d=4$ or $d\ge 6$. We show that the minimal number of required pieces is less than $20d$ if $d \ge 10$.

\end{abstract}
\section{Introduction}

The history of this problem goes back to $1949$, when Van der Waerden posed an exercise in Elemente der Mathematik. The question was whether the disk can be decomposed into 2 disjoint congruent pieces. Different elementary proofs show that it is not possible.

Maybe the simplest one is the following: If there exists such a decomposition, then there exists an isometry connecting the two pieces. We prove that this isometry must be a linear transformation. Let $A \cup B =D$ be a decomposition of the unit ball and $\phi$ be an isometry with $\phi(A)=B$. The 1 dimensional Hausdorff measure of the boundary of the disk, $\mathcal{H}_1(\partial D)$ is $2 \pi$. The outer Hausdorff measure of the intersection of the boundary of the disc with $A$ or $B$ is at least $\pi$. We may assume that this holds for $A$. On the other hand, there is no arc of radius 1 contained in the interior of the disk which has at least $\pi$ measure. Therefore $\phi(A \cup \partial D) \subset \partial D$. Then the origin stays in place.

This motivates the question whether the $d$ dimensional ball can be decomposed into finitely many congruent pieces. Clearly, it is enough to decide the question for the unit ball $B_d$. For a cardinal number $m$, we say that a set $K$ is {\bf $m$-divisible} (with respect to $G$) if $K$ can be decomposed into $m$ congruent (with respect to $G$) and disjoint pieces.   Wagon \cite{W1} proved in 1984 that the $d$ dimensional ball is not $m$-divisible for $2 \le m \le d$. This was the only well-known lower bound for the number of pieces. In 2012 the authors showed (in an unpublished paper) that the disk is not $3$-divisible.

In 2007, Richter \cite{R} showed that a typical convex body $D$ is not $m$-divisible for any finite $m$. Every decomposition can be described by a set $A$ and a set of isometries $\phi_0=id, \phi_1 \stb \phi_n$, where $D= \coprod_{i=0}^{n} \phi_i(A)$.
He proved that if $\mathcal{H}_{d-1}(\phi_i^{-1}(\partial D) \cap \phi_j^{-1}(\partial D))=0$ for every $i \ne j$, then $D$ cannot be decomposed by these isometries. This guarantees that every element of a residual subset of the space of convex bodies (endowed with Hausdorff metric) is not $m$-divisible for any $m\ge 2$. However, for every $d$, the $d$ dimensional ball $B_d$ is not in this class, see \cite{R}. In 2010, Laczkovich and the first author proved \cite{KL} that the 3 dimensional ball is $m$-divisible for any $m\ge 22$.

In this paper we prove that the $d=2k$ dimensional ball can be decomposed into finitely many congruent pieces:
\begin{theorem}\label{t12}
The $2s$ dimensional ball (either open or closed) is $m$-divisible
for every $m \ge 4(2s+1)+2$ if $s \ge 2$ and $s \ne 3$.
\end{theorem}
The original proof was formulated for the four dimensional unit ball. The construction of the proof is a natural generalization of it for higher dimensional cases. As a special case of Theorem \ref{t12} we get:

\begin{theorem}\label{t11}
The $4$ dimensional ball (either open or closed) is $m$-divisible
for every $m \ge 22$.
\end{theorem}

 Using Theorem \ref{t12} and the fact that the 3 dimensional ball can be decomposed into finitely many pieces
  (see \cite{KL}), we prove the following:
\begin{theorem}\label{t2}
The $d$-dimensional ball $B_d$ can be decomposed into finitely many pieces for $d \ge 6$ and $d=3,4$.
\end{theorem}

Furthermore, we show that the minimal number of pieces in our construction grows linearly with the dimension. According to \cite{W1} this is the best in the sense that there is a linear lower bound $d+1$ for the number of pieces which is needed for a decomposition.

\begin{theorem} \label{t3}
Let $d\ge 10$ and $\tau(B_d)$ denote the minimal number of required pieces for a decomposition of the ball $B_d$.
Then $$ d< \tau(B_d) < 20\cdot d.$$
\end{theorem}

This result improves the upper bound given by the construction for the $3k$-dimensional ball in \cite{KL}.

Our paper is organized as follows. In Section \ref{notation} we introduce the notation that we will use throughout the paper.
In Section \ref{s3} we collect facts about a subgroup of the $4$ dimensional special orthogonal group. In Section \ref{independence} we define a rational parametrization of special orthogonal matrices in dimensional $d$. Section \ref{4k} is devoted to the main lemma of the paper giving sufficient properties for the existence of decomposition of infinite graphs. In Section \ref{s6} and \ref{s7} we apply it for graphs defined by isometries. In Section \ref{s6} we complete the proof of Theorem \ref{t11} and Theorem \ref{t12}. Finally, in Section \ref{s7} we handle the odd dimensional cases to prove Theorem \ref{t2} and we collect all the information given in the paper on the number of required pieces for a decompositions to prove Theorem \ref{t3}. In Section \ref{vegsos8} we summarize the results and open questions on the decomposition of balls.

\section{Notation}\label{notation}
For a possible directed graph $\Gamma$ we denote by $V(\Gamma)$ and $E(\Gamma)$ the set of vertices and edges, respectively. If there is an edge $e$ from $U$ to $V$, then we say that $U$ is the tail and $V$ is the head of $e$ and we denote them by $T(e)$ and $H(e)$, respectively.
We call a sequence of vertices $V_1, V_2 \stb V_n$ a path if for every $1 \le i \le n-1$ there is an edge from $V_i$ to $V_{i+1}$ and a path $V_1, V_2 \stb V_n$ is a cycle if $V_i \ne V_j$ if $1 \le i <j \le n-1$ and $V_1=V_n$.
We denote by $(P,Q)$ a path from $P$ to $Q$. We will also use this notation for graphs, where there are more than one paths connecting $P$ and $Q$ if it is clear which path we consider.

We denote by $e$ the identity element of a group. Let $G$ be a group generated by the elements of the set $S=\{ w_{\alpha} \mid \mbox{ } \alpha\in I \}$, where $S=S^{-1}$.  Every element $W$ of the group $G$ can be written (not necessarily uniquely) as a word of the generators so $W$ is of the form $w_1 w_2 \ldots w_n$, where $w_i \in S$ for every $1 \le i \le n$. As a later terminology, we say that the word $W$ starts with $w_1$ and ends with $w_n$. Moreover the i'th letter $w_i$ of the word $W$ will be denoted by $W[i]$ and we use the notation $W[-1]$ for the last letter of $W$. If $W$ is the empty word, then let $W[i]=e$. We denote by $\lg(W)$ the length $n$ of the reduced word $W$.
However, we will use linear transformations of $\rr^d$ as the letters of a word $W$ and we use the convention that linear transformations acts from the left on the elements of $\rr^{d}$.
We also say that a word $W=w_1 w_2 \stb w_k$ has a shorter conjugate if $W[1]^{-1}=W[-1]$.

The special orthogonal group $SO(n,\rr)$ will be shortly denoted by $SO(n)$ and we denote by $Iso(n)$ the isometry group of the $n$ dimensional Euclidean space. In this paper, by $m$-divisibility we mean $m$-divisibility with respect to $Iso(n)$.

Let $p(x) =a_n x^n + \dots +a_0$ be a polynomial. Let $\deg(p)$ denote the degree $n$ of the polynomial $p$ and we denote by $LC(p)$ the leading coefficient $a_n$ of $p$.

\section{Lemmas on a subgroup of $SO(4)$}\label{s3}
In this section, for sake of completeness, we prove more than it would be necessary to prove Theorem \ref{t11}.
\begin{lemma} \label{l40}
Let $A$ and $B$ be the rotations in $SO(4)$ given by the matrices
\begin{equation*} \begin{split} A&=\left(
\begin{matrix}
\cos \theta & -\sin \theta & 0 & 0\\
\sin \theta & \phantom{-} \cos \theta & 0 & 0 \\
0 & \phantom{-} 0 & \cos \theta & -\sin \theta \\
0 & 0 & \sin \theta & \phantom{-} \cos \theta \\
\end{matrix}
\right)
 \text{ and} \\
B&=\left(
\begin{matrix}
\cos \theta & 0 & 0 & -\sin \theta\\
0& \cos \theta & -\sin \theta & 0 \\
0& \sin \theta & \phantom{-} \cos \theta & 0 \\
\sin \theta & 0 & 0 & \phantom{-} \cos \theta\\
\end{matrix}
\right) ,
\end{split} \end{equation*}
respectively, where $\cos \theta$ is transcendental. We denote by $K$ the group generated by $A$ and $B$.
Then every element $U\ne 1 \in K$ has exactly one fix point, which is the origin.
\end{lemma}

\proof
The proof can be found in \cite[Theorem 6.3]{W2}.
$\square$

\begin{obs}\label{szabadhossz}
It is easy to see from Lemma \ref{l40} that $K$ is a free group so every element of $K$ can be written uniquely as the product of the matrices $A, A^{-1}, B, B^{-1}$. This gives that the length of $M =A^{m_1} B^{n_1} \cdots A^{m_l} B^{n_l} \in K$ is $\lg(M) =\sum_{i=1}^{l} |m_i|+|n_i|$ if $M$ is defined by a reduced word.
\end{obs}
\begin{defin}
\begin{enumerate}
\item We define the set $$ \mathcal{M} =\left\{ \begin{pmatrix}
a & -b & -c & -d\\
b & a & -d & c \\
c & d & a & -b \\
d & -c& b & a\\
 \end{pmatrix} : a, b, c, d \in \rr
 \right\}$$
It is easy to verify that $\mathcal{M}$ is an algebra over $\rr$.
\item Let
$$ \mathcal{M}_1 =\left\{ M \in \mathcal{M} : \mbox{ } det(M) =1  \right\} .$$
Clearly, $\mathcal{M}_1$ is a subgroup of the orthogonal group $O(n)$.
\item
 Similarly, let $\mathcal{M}(\theta)$ denote the set of matrices of the form
 \begin{equation*}\begin{split} \left( \begin{matrix}
p(\cos\theta) & -\sin\theta ~q(\cos\theta) & -r(\cos\theta  ) & -\sin\theta ~s(\cos\theta)\\
\sin\theta ~q(\cos\theta) & p(\cos\theta) & -\sin\theta ~s(\cos\theta) & r(\cos\theta  ) \\
r(\cos\theta ) & \sin\theta ~s(\cos \theta) & p(\cos \theta) & -\sin \theta ~q(\cos\theta) \\
\sin\theta ~s(\cos \theta) & -r( \cos\theta ) & \sin\theta ~q(\cos\theta) & p(\cos\theta)\\
 \end{matrix}   \right) \mbox{,}
 \end{split}
  \end{equation*}
 where $p, q, r, s \in \mathbb{Q}[x]$.
  Such an element of $\mathcal{M}(\theta)$ is determined by the polynomials $p, q, r, s$ and will be denoted by $M_{\theta}(p,q,r,s)$.
\item
Let $\mathcal{M}_1(\theta)=\{M\in \mathcal{M}(\theta): \det{M}=1\}$.
 \end{enumerate}
 \end{defin}

 \begin{obs}\label{l41}
Let $U$ be an element of $K$, where $K$ is defined in Lemma \ref{l40}.
 Then $U \in \mathcal{M}_1(\theta)$.
 \end{obs}

For further results we need to describe the degree and leading coefficient of the polynomials $p,q,r,s$ for $M_{\theta}(p,q,r,s) \in \mathcal{M}(\theta)$.

\begin{defin}
\begin{enumerate}
\item\label{d4.2}
For a pair of polynomials $p_1$, $p_2$ we write $p_1(\cos (\theta)) \doteq p_2(\cos (\theta ) )$ and $\sin(\theta) p_1(\cos (\theta)) \doteq \sin(\theta) p_2(\cos (\theta ) )$ if $\deg(p_1) = \deg(p_2)$ and $LC(p_1) =LC(p_2)$.
\item
For a pair of matrices $M_1, M_2 \in \mathcal{M}(\theta)$ we write $M_1 \doteq M_2$ if and only if
$(M_1)_{i,j}\doteq(M_2)_{i,j}$ for every $i,j\in\{1\stb 4\}$.

\end{enumerate}
\end{defin}
We define the degree of a matrix in $M(\theta)$.
\begin{defin}\label{matrixdegree}
Let $M=M_{\theta}(p,q,r,s)$. We denote by $\deg(M)$ the maximum of $\deg(p),\deg(q)+1,\deg(r),\deg(s)+1$.
\end{defin}
It is easy to see that if $M,N \in \mathcal{M}(\theta)$ and $M \doteq N$, then $\deg(M)= \deg(N)$.
\begin{obs}\label{dot}
\begin{enumerate}[(a)]
\item\label{dot1}
It is easy to see that for $p(x) =a_n x^n + a_{n-1} x^{n-1} \stp a_0$ we have $p(\cos \theta) \doteq a_n (\cos \theta)^n$. We also have $\sin \theta p(\cos \theta) \doteq a_n \sin \theta (\cos \theta)^n$.
\item\label{dot2}
Let $p_1$, $p_2$, $q_1$ and $q_2$ polynomials in $\mathbb{Z}[x]$.
If $p_1 \doteq q_1$ and $p_2 \doteq q_2$, then $p_1 p_2 \doteq q_1 q_2$.
\item\label{dot3}
 Let us assume again that $p_1 \doteq q_1$ and $p_2 \doteq q_2$. If $\max \left\{ \deg(p_1), \deg(p_2) \right\}=\deg(p_1+ p_2)$, then $  \max \left\{ \deg(q_1) , \deg(q_2) \right\}=\deg(q_1 + q_2)$ and $p_1 +p_2 \doteq q_1 +q_2$.
\item\label{dot4}
If $\deg(p_1) > \deg(p_2)$, then $p_1+p_2 \doteq p_1$.
\end{enumerate}
\end{obs}

\begin{lemma}\label{l401}
Let $U \in K$ be of the form  $A^{m_1}B^{n_1} \cdots A^{m_t}B^{n_t}$, where $A$ and $B$ are given in Lemma \ref{l40}.
Let $\sigma$ denote the length of $U$.
\begin{enumerate}[(a)]
\item\label{a} (Case $U=A^{m_1}B^{n_1} \cdots A^{m_t}B^{n_t}$)\\
If $m_i ,n_i$ are nonzero integers for $1 \le i \le t$,
then $U \in \mathcal{M}(\theta)$, where $\deg(p)=\deg(r)=\deg(q)+1=\deg(s)+1=\sigma$. We also have $|LC(p)|=|LC(q)|=|LC(r)|=|LC(s)|=2^{\sigma -t-1}$.
\item\label{bb} (Case $U=A^{m_1}B^{n_1} \cdots A^{m_t}$)\\
If $m_i ,n_i, m_t$ are nonzero integers for $1 \le i \le t-1$  and $n_t = 0$, then one of the following two cases holds:
\begin{enumerate}[label=\roman*]
\item  $\deg(p)=\deg(q)+1= \sigma$ with
$|LC(p)|=|LC(q)|=2^{\sigma -t-2}$ \\and   $\max(\deg(r),\deg(s)+1)< \sigma$.
\item $\deg(r)=\deg(s)+1= \sigma$ with $|LC(r)|=|LC(s)|=2^{\sigma -t-2}$ \\and   $\max(\deg(p),\deg(q)+1)< \sigma$.
\end{enumerate}
\item\label{c} (Case $U=B^{n_1} \cdots A^{m_t}B^{n_t}$, similarly)\\
If $n_1 ,m_i, n_i$ are nonzero integers for every
$2 \le i \le t$  and $m_1 = 0$, then one of the following two cases holds:
\begin{enumerate}[label=\roman*]
\item  $\deg(p)=\deg(s)+1= \sigma$ with
$|LC(P)|=|LC(S)|=2^{\sigma -t-2}$ \\and $\max(\deg(r),\deg(q)+1)< \sigma$.
\item $\deg(q)+1=\deg(r)= \sigma$ with $|LC(q)|=|LC(r)|=2^{\sigma -t-2}$ \\
and   $\max(\deg(p),\deg(s)+1)< \sigma$.

\end{enumerate}
\end{enumerate}
\end{lemma}

\proof
\begin{enumerate}[(a)]
\item
We claim that
\begin{equation}\label{indula}
U \doteq 2^{\sigma -t- 1} \cos^{\sigma -1} \theta \left(
\begin{matrix}
\xi \cos\theta & -\mu \sin \theta & - \zeta \cos\theta & - \nu \sin \theta \\
 \mu \sin \theta  & \xi \cos\theta & -\nu \sin\theta & \zeta \cos\theta  \\
\zeta \cos\theta  & \nu \sin\theta & \xi \cos\theta & -\mu \sin \theta \\
\nu \sin\theta & -\zeta \cos\theta  & \mu \sin \theta & \xi \cos\theta \\
\end{matrix} \right) \end{equation}
for some $\xi, \mu, \zeta, \nu= \pm 1$ with $\mu \nu = \zeta \xi$.
The proof of this fact can be found in Wagon \cite[page 55]{W2}.
\item We write $U=U'A^{m_t}$, where $U'= A^{m_1}B^{n_1} \cdots A^{m_{t-1} }B^{n_{t-1} }$. Equation \eqref{indula} shows that \begin{equation*}
U' \doteq 2^{\sigma' -t- 2} \cos^{\sigma' -1} \theta \left(
\begin{matrix}
\xi' \cos\theta & -\mu' \sin \theta & - \zeta' \cos\theta & - \nu' \sin\theta \\
 \mu' \sin \theta  & \xi' \cos\theta & -\nu' \sin\theta & \zeta' \cos\theta  \\
\zeta' \cos\theta  & \nu' \sin\theta & \xi' \cos\theta & -\mu' \sin \theta \\
\nu' \sin\theta & -\zeta' \cos\theta  & \mu' \sin \theta & \xi' \cos\theta \\
\end{matrix} \right), \end{equation*}

where $\sigma'= |m_1|+|n_1|+\dots+|m_{t-1}|+|n_{t-1}|$.
Using the fact that $\mu' \nu' = \zeta' \xi'$ and $|\mu'| =|\nu'| = |\zeta'|= |\xi'|$ we get that exactly one of the two sums $\xi'+ \mu'$ and $\zeta'-\nu'$ is 0 and the absolute value of the other is $2$.

It is easy to show that
\begin{equation*}
\begin{split}
A^{m_t}&= \left(
\begin{matrix}
\cos (m_t \theta) & -\sin (m_t \theta) & 0 & 0\\
\sin (m_t \theta) & \phantom{-} \cos (m_t \theta) & 0 & 0 \\
0 & \phantom{-} 0 & \cos (m_t \theta) & -\sin (m_t \theta) \\
0 & 0 & \sin (m_t \theta) & \phantom{-} \cos (m_t \theta) \\
\end{matrix}
\right)  \\ &\doteq 2^{m_t-1} \cos^{m_t-1} \theta ~\left(
\begin{matrix}
\cos \theta & -\sin \theta & 0 & 0\\
\sin \theta & \phantom{-} \cos \theta & 0 & 0 \\
0 & \phantom{-} 0 & \cos \theta & -\sin \theta \\
0 & 0 & \sin \theta & \phantom{-} \cos \theta \\
\end{matrix}
\right),
\end{split}
\end{equation*}

using well known facts about Chebyshev polynomials. \\ We define the matrix $M$ and $N$ by $U'=(2^{\sigma' -t- 2} \cos^{\sigma' -1} \theta)\cdot M$ and $A^{m_t} =(2^{m_t-1} \cos^{m_t-1} \theta)\cdot N$.
Using the fact $\sin^2 \theta =1-\cos^2 \theta$, we get that the first row of $M\cdot N$, which is denoted by $(M\cdot N)_{1\cdot}$, is the following
\begin{equation*} \begin{split} (M\cdot N)_{1\cdot} \doteq & \left( (\xi'+\mu')cos^2\theta,-(\xi'+\mu')sin\theta\cdot cos\theta, \right. \\ & \left.-(\zeta'-\nu')cos^2\theta,(\zeta'-\nu')sin\theta\cdot cos\theta \right) \mbox{.} \end{split} \end{equation*}
Thus either the first two or the second two coordinates vanishes. Easy calculation shows that in the other two coordinates of $U' A^{m_t} $ have degree  $(\sigma'-1)+(m_t-1)+2=\sigma$ and $\sigma-1$, respectively. The absolute value of the leading coefficients are the same $(2^{m_t-1}2^{\sigma' -t- 2}\cdot 2 =2^{\sigma-2}$).
\item Similar calculation shows the statement.
\end{enumerate}
\qed
\begin{lemma}\label{foatlo}
Let $U$ and $\sigma$ be as in Lemma \ref{l401} \ref{bb}.
We claim that $m_1 m_t >0$ if and only if $$\deg(p) = \deg (q)+1 \mbox{ and } \max (\deg(r), \deg(s)+1) < \deg(p)\mbox{.}$$
\end{lemma}

\proof
If $M=M_\theta(p,q,r,s) \in \mathcal{M}(\theta)$, then $tr(M) =4 p(\cos(\theta) )$.
Conjugating by $A^{m_t}$ we get
$$tr(M)=tr \left( A^{m_1}B^{n_1} \cdots B^{m_{t-1}}A^{m_t}\right) = tr(A^{m_1+m_t}B^{n_1} \cdots B^{m_{t-1}}) \mbox{.}$$
Clearly, the sum of the absolute value of the exponents $\sigma' = |m_1+m_t| +|n_1| \stp |m_{t-1}|$ is smaller than $\sigma-1$ if $m_1 m_t<0$ and $\sigma' = \sigma$ if $m_1 m_t>0$. By Lemma \ref{l401} \ref{a} we have $\deg(p) = \sigma' < \sigma-1$ if $m_1 m_t<0$ and $\deg(p) =\sigma'= \sigma$ if $m_1 m_t>0$. Finally, one can identify the two cases of Lemma \ref{l401} \ref{bb}, finishes the proof.
\qed

\begin{remark}
\begin{enumerate}
\item
The analogue statement is true for $U$ and $\sigma$ in Lemma \ref{l401} \ref{c} and for $n_1, n_t$ instead of $m_1, m_t$. Therefore, for a matrix $$M\in\mathcal{M}(\theta) $$
$\deg(M)$, which is defined in Definition \ref{matrixdegree}, is taken in the diagonal if and only if $M$ does not have a shorter conjugate.
\item
Now we can easily calculate the degree of the polynomials in the main diagonal of the word $U$ which equals to
\[ \min \left\{ lg(U') \mid U' \in K, \mbox{ } U \mbox{ and } U' \mbox{are conjugate} \right\}   \mbox{.} \]
\end{enumerate}
\end{remark}
It is easy to see that every element $M$ of the group generated by $A$ and $B$ we have $\lg(M)=\sigma =\deg(M)$.

\begin{lemma}\label{inverz}
Suppose that $M \in \mathcal{M}$ and $tr(M)=4p\ne 1$. Then $M+M^{T} =2pI$ is a scalar matrix and $(I-M)^{-1} = \frac{1}{2-2p} (I-M^{T})$.
\end{lemma}
\proof
Clearly, $M+M^{T} =2pI$ and $$(I-M)(I-M^{T}) = I-M-M^{T} + MM^{T} = I-2p I+I=(2-2p)I$$ since $M$ is an orthogonal matrix.
\qed
\\
Technically, we need the following as well:
\begin{lemma}\label{inverzaltalaban}
Suppose that for $M \in \rr^{d \times d}$, the matrix $I-M$ is invertible, then the entries of $(I-M)^{-1}$ are rational functions of the entries of $M$.
\end{lemma}
\proof
Obvious, using Cramer's rule.
\qed
\section{Algebraic independence}\label{independence}

It was proved in \cite[p. 5-6.]{KL} that there exists a rational parametrization $\alpha_d : \Omega \rightarrow SO(d)$, where $\Omega$ is an open subset of $\mathbb{R}^{d'\cdot d}$ and $\alpha_d$ is surjective, where $d'=d$ if $d$ is even and $d'=d-1$ if $d$ is odd.
Indeed, every element of $SO(d)$ can be written as the product of at most $d'$ reflections given by the vectors $ v_i=(x_{di+1}\stb x_{di+d} )$ for $i=0, 1 \stb d'-1$. For every $w\in \rr^d$  the matrices $$R_w=I-\frac{ww^T}{|w|^2}$$ gives a parametrization of the reflection in a hyperplane perpendicular to $w$. The entries of the matrix $R_w$ are rational functions of the coordinates of $w$, where the denominator of the functions does not vanish for any $w\ne {\bf 0}$. Let $v=(v_0,v_1 \stb v_{d'-1})$, which is the concatenation of the vectors $v_i \in \rr^d$.
Hence the entries of the matrix $\alpha(v)$ are rational functions of $x_1, x_2 \stb x_{d'd}$ with integer coefficients. The denominator of $(\alpha_d (v))_{i,j}$ does not vanish on $\Omega$ as a rational function.

Now, we fix a rational parametrization $\alpha_d$ of $SO(d) $. If $v\in \Om
\su \rr^{d' d}$,
then we shall denote by $O_v$ the image of the parametrization;
both as a matrix and as a linear transformation of $\rr ^d $.
Then $v \mapsto O_v$
is a surjection from $\Om$ onto $SO(d)$, and every entry of
the matrix of $O_v$ is a
rational function with integer coefficients of the coordinates of $v$.

\begin{defin}\label{independentdef}
We say that $M_1, M_2 \stb M_m \in SO(d)$ are independent, if there exist $v_1, v_2  \stb v_m \in \Omega$ such that $\alpha_d(v_i)=M_i$ and the coordinates of $v_i$ are algebraically independent over $\mathbb{Q}$.
We will also say that a vector $t \in \rr^{d}$ and the matrices $M_1, M_2 \stb M_m \in SO(d)$ form an independent system if the coordinates of $t$ and the coordinates of $v_1 \stb v_m$ are algebraically independent over $\mathbb{Q}$.
\end{defin}

\begin{lemma}\label{ftln}
Let $p$ be a polynomial on $k$ variables. Let $M_1 \stb M_k$ be independent elements of $SO(d)$ with $p(M_1 \stb M_k)=0$. Then  $p(N_1\stb N_k)=0$ for all $N_1 \stb N_k \in SO(d)$.
\end{lemma}
\proof
Every entry of the matrix equation is a polynomial expression of the parameters. Since they were chosen algebraically independently, the equation holds if and only if it is trivial. This means that it holds for any substitution of the parameters. The fact that $\alpha_d$ is surjective finishes the proof of Lemma \ref{ftln}.
\qed

Similar argument shows the following.
\begin{lemma}\label{ftlnremark}  Let $q(x_1\stb x_{N+d})$ be a rational function where $N=d^2\cdot k$ and $M_1 \stb M_k  \in SO(d)$  and $t \in \rr^d$ be an independent system. Let us suppose that $q(m_1\stb m_N, t_1\stb t_d)=0$ where $(m_i)$ is an enumeration of the entries of the matrices $M_1\stb M_k$. Then $q(n_1\stb n_N,s_1\stb s_d)=0$ holds for the same enumeration of the entries $(n_i)$ of arbitrary matrices $N_1 \stb N_k \in SO(d)$ and arbitrary $s=(s_1\stb s_d)\in \rr^d$, where the denominator of $q$ does not vanish.
\end{lemma}

We usually use this fact contrary, we show that there exists a substitution which is non-trivial, therefore it is non-trivial for any algebraically independent substitution.

\begin{lemma}\label{nincsfixpont}
Let $\alpha_d$ be a rational parametrization of the $d$ dimensional special orthogonal linear transformations, where $d$ is even. If $O_1, O_2 \stb O_k$ are independent orthogonal transformations and $U$ is not an empty word, then $\hat{U}=U(O_1, O_2 \stb O_k)$ does not have a nonzero fix point.
\end{lemma}
\proof
The characteristic polynomial $p(y)=\det(Iy-U)$ of the orthogonal transformation $\hat{U}=U(O_1, O_2 \stb O_k)$ can be considered as a rational function with integer coefficients of the variables $y, x_1,x_2 \stb x_{k\cdot d^2}$.
 Let us assume indirectly that $\hat{U}$ has a nonzero fixpoint, thus $p$ vanishes at $y=1$. By the algebraic independence of the parameters we get that $p$ vanishes at $y=1$ for any substitution to the variables $x_1, x_2 \stb x_{k\cdot d^2}$. This shows that $1$ is the eigenvalue of every element of the form $U(M_1, M_2 \stb M_k)$, where $M_i$ are orthogonal transformations, which clearly contradicts Lemma \ref{l40} if $d=4$. Moreover, free subgroup of the orthogonal group consisting of fixed point free elements (except the identity) was given in \cite{Dekker,DS} for every $d$ dimensional orthogonal groups where $d$ is even and $d\ge 4$, finishing the proof of Lemma \ref{nincsfixpont}.
\qed
\section{Decomposition in $\mathbb{R}^{2s}$}\label{4k}
Let $X$ be a set, and let $f_1 \stb f_n$ be maps from subsets of $X$ into $X.$
Our aim is to find a sufficient condition for the existence of a decomposition
$X=A_0 \cup A_1 \cup \ldots \cup A_n$
such that $f_i (A_0)=A_i$ for every $i=1\stb n.$

Suppose that for $i=1, 2 \stb n$ the function $f_i$ is defined on $D_i \su X$ $(i=1\stb n),$
and put $D=\bigcap_{i=1}^n D_i .$
We say that the point $x$ is a core point, if $x\in D,$
and the points $x,f_1 (x)\stb f_n (x)$ are distinct.
By the image of a point $x$ we mean the multiset $\mathcal{I}_x=\{ f_1 (x) \stb f_n (x)\}$. The multiset $\mathcal{I}_x$ is a set if $x$ is a core point.

For a set $\mathcal{F}=\left\{f_1 \stb f_n \right\}$ we define a graph $\Ga_{\mathcal{F}}$ on the set $X$ as follows.
We connect the
distinct points $x,y\in X$ by an edge if there is an $i\in \{ 1\stb n\}$ such
that $f_i (x)=y$. Then $\Ga_{\mathcal{F}}$ will be called the graph
generated by the functions $f_1 \stb f_n .$

\begin{lemma}\label{l1}
Let $X,$ $f_1  \stb f_n ,$ $D,$ and $\Ga_{\mathcal{F}}$ be as above,
and suppose that the graph $\Ga_{\mathcal{F}}$
has the property that
\begin{equation}\label{coin}
 \begin{split}
 &\textrm{whenever two cycles } \mathcal{C}_1 \textrm{ and } \mathcal{C}_2 \textrm{ in }
 \Ga_{\mathcal{F}} \textrm{ share a common edge}, \\
&\textrm{ then the sets of vertices of } \mathcal{C}_1 \textrm{ and } \mathcal{C}_2
\textrm{ coincide}.
\end{split}
\end{equation}
Suppose further that there is a point $x_0 \in X$
satisfying the following conditions.
\begin{enumerate}[(a)]
\item\label{fofeltetel1}
$x_0$ is in the image of at least one core point;
\item\label{fofeltetel2}
every $x\in X\se \{ x_0 \} $
is in the image of at least three core points.
\end{enumerate}
Then there is a decomposition $X=A_0 \cup A_1 \cup \ldots \cup A_n$
such that $A_0 \su D,$ and
$f_i (A_0)=A_i$ for every $i=1\stb n.$
\end{lemma}
 \proof The proof is based on the axiom of choice and can be found in \cite{KL}.
 \qed

\begin{lemma}\label{uj1}
If a connected component $\Gamma'$ of $\Gamma$ contains two different cycles sharing at least a common edge, then $\Ga'$ contains two cycles $\mathcal{C}_1 = P_1, P_2 \stb P_m $ and $\mathcal{C}_2= Q_1, Q_2 \stb Q_n$ such that for some $1 < k < \min \{n,m \}$ we have $P_i =Q_i$ for $i=1 \stb k$ and  $\{ P_1, P_2 \stb P_m \} \cap \{ Q_1, Q_2 \stb Q_n \} = \{ P_1, P_2 \stb P_k \}$.
\end{lemma}
\proof
We may assume that $P_1 =Q_1$ is one of the endpoints of a common edge such that $P_2 \ne Q_2$. Then there exists a minimal integer $b$ such that $Q_b =P_a$ for some $1 <a < m$. Since $P_1,P_2 \stb P_m$ are different points, the cycles $Q_1, \ldots ,Q_b, P_{a-1}, \ldots , P_1$ and $Q_1, \ldots ,Q_b, P_{a+1}, \ldots , P_m$ have a common path and share only the points $Q_1, Q_2, \ldots , Q_b$.
\qed

\begin{remark}\label{minsys}
Essentially, this means that we can find two points $P$ and $Q$ such that between these points there are three paths which have no other common points.
\end{remark}

\begin{theorem}\label{foprop}
Let us assume that $t \in \rr^{d}$ and $\alpha_{v_0}= O_0, \alpha_{v_1}=O_1 \stb \alpha_{v_m}=O_m$ in $SO(d)$ form an independent system, where $d=2s \ge 4$ and $d\ne 6$. Let $F(x) = O_0 x +t$. Then $\Ga_{\mathcal{F} }$ has the property \eqref{coin}, where $\mathcal{F} = \left\{ F, O_1 \stb O_m \right\}$.
\end{theorem}

\proof
Let us assume indirectly that there exists a connected component of $\Ga'$ which the contains cycles $\mathcal{C}_1$, $\mathcal{C}_2$ and the two cycles share at least one edge. Using Lemma \ref{uj1} we may assume that the two cycles share a common path. Thus $\Ga'$ contains a subgraph $\Delta=(V(\Delta), E(\Delta))$:

\begin{center}
\includegraphics[width=8cm]{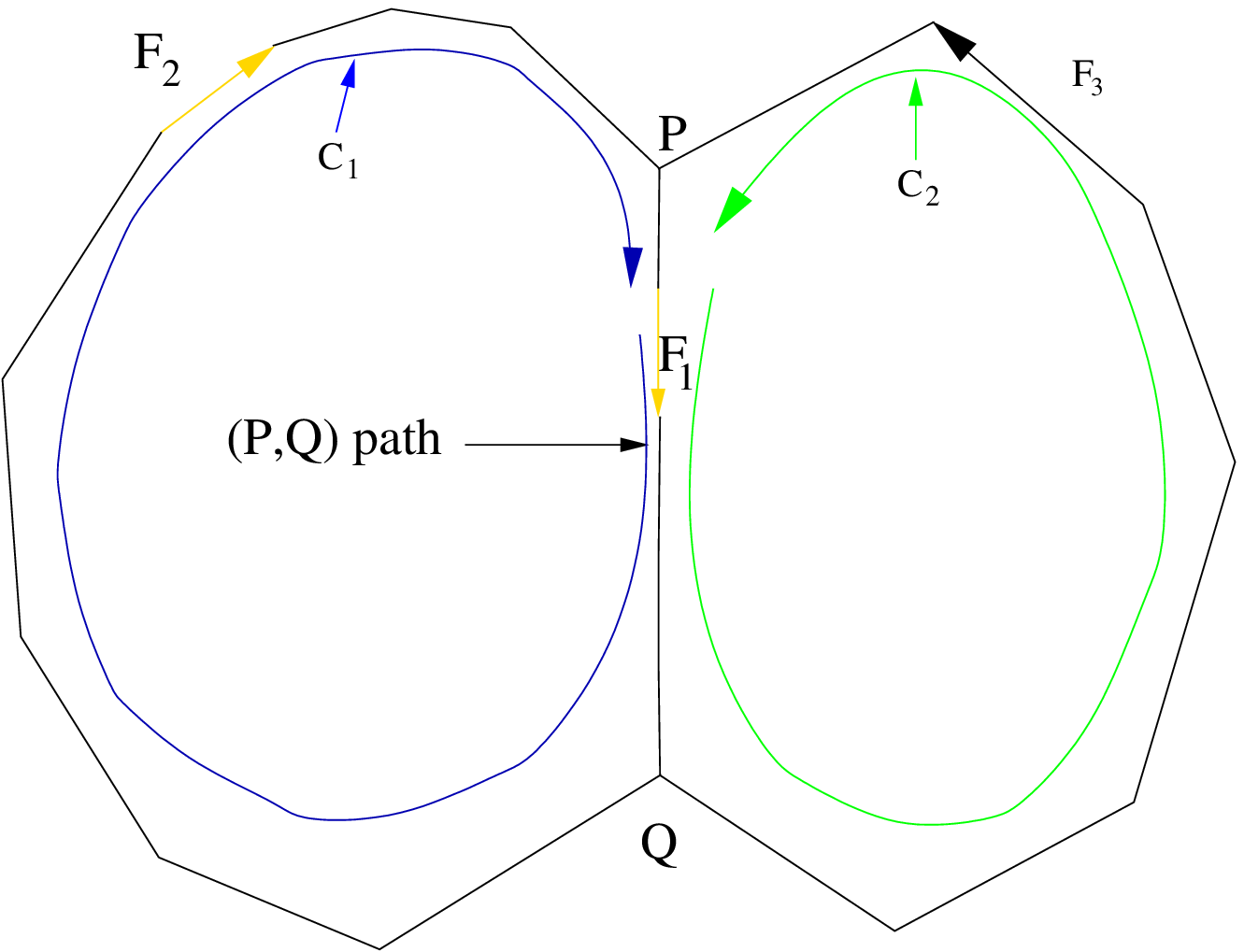}
\end{center}

Let us denote by $P$ and $Q$ the endpoints of the common paths and denote by $(P,Q)$ path the common path as in Figure 1.
For each edge of the graph we can naturally assign a letter $O_i$ or $F$.
\begin{remark}\label{F akoron}
$O_i$ are independent orthogonal transformations and $F$ is the only isometry involving translation, therefore by Lemma \ref{nincsfixpont} there must be a letter $F$ or $F^{-1}$ in every cycle.
Thus we may assume that at least two of the three paths between $P$ and $Q$ contain the letter $F^{\pm 1}$.
\end{remark}
Using the previous remark we may assume that the $(P,Q)$ path contains an $F$ or an $F^{-1}$. We denote the closest $F^{\pm 1}$ to $P$ on the path $(P,Q)$ by $F_1$.

We choose a starting point $S$ from which we start going around the cycles $\mathcal{C}_1, \mathcal{C}_2$ (as in Figure 1) and then the two cycles naturally determine two words $W_1$ and $W_2$, respectively.
 According to Remark \ref{F akoron}, there is another $F^{\pm 1}$ in $W_1$, which as an edge is not contained in $E(\mathcal{C}_2)$. Let us denote the first $F^{\pm 1}$ in $W_1$ by $F_2$. Similarly to $F_1$ and $F_2$ one can define $F_3$ to be the edge corresponding to the last $F^{\pm 1}$ on the cycle $\mathcal{C}_2$. Note that $F_3$ might be equal to $F_1$ and it might also happen that $F_1 \ne F_3$ but $F_3$ is on the $(P,Q)$ path. We consider the edge corresponding to $F_2$ and $F_3$ as a directed edge which has the same direction as the cycle $\mathcal{C}_1$ and $\mathcal{C}_2$, respectively.

The starting point $S$ of the two cycles can be identified with an element of $x \in \rr^{d}$ which satisfies
\begin{equation}\label{xS} W_1(x)=W_2(x)=x \mbox{.} \end{equation}

Every direction-preserving isometry of $\rr^d$ can be written as $W(x)=U(x)+b$ for some $U \in SO(d)$ and $b \in \rr^d$.
Using equation \eqref{xS} we get that there are $U_1, U_2\in SO(d)$ and $b_1, b_2 \in \rr^d$ such that
\begin{equation}\label{ketkor}
  W_1(x)=U_1x+b_1=x  \mbox{ and } W_2(x)=U_2x+b_2 =x \mbox{.}
\end{equation}
Let $H$ denote the group generated by $O_0, O_1 \stb O_m$. Since the edges of $\Delta$ are labelled by $F^{\pm 1}$ and $O_i^{\pm1}$ we have $U_1$ and $U_2$ are in $H$. Thus
$$ (I-U_1)x=b_1 \mbox{ and } (I-U_2)x=b_2 \mbox{.} $$

Since $\mathcal{C}_1$ and $\mathcal{C}_2$ are cycles, $U_1$ and $U_2$ are nonempty reduced words of the generators $O_0, \stb O_m$.
By Lemma \ref{nincsfixpont}, the orthogonal transformation $U_1$ and $U_2$ do not have a fix point thus $I-U_i$ are invertible for $i=1,2$ and hence
\begin{equation}\label{e2x}
(I-U_1)^{-1}b_1= (I-U_2)^{-1}b_2=x \mbox{.}
\end{equation}
One can easily verify that $(I-U_1)^{-1}b_1= (I-U_2)^{-1}b_2$ is equivalent to the fact the words $W_1$ and $W_2$ have a common fix point, which was formulated in equations \eqref{xS} and \eqref{ketkor}.

We write
 \[ W_i(x) =S_{i,1} F^{\alpha_{i,1} } S_{i,2}  F^{\alpha_{i,2} } \dots S_{i,n_i} F^{\alpha_{i,n_i} }S_i^{*}(x) \mbox{,} \]
where $S_i^{*}$ and $S_{i,j}$ are elements of the group $H'=\langle O_1 \stb O_m \rangle$ and $\alpha_{i,j}$ is $1$ or $-1$ for every $j=1 \stb n$ and $i=1,2$.
In this case for $W_i$ is of the following form for $i=1,2$:
\begin{equation}\label{hosszu}
\begin{split} &W_i=S_{i,1} {O_0}^{\alpha_{i,1} } S_{i,2} {O_0}^{\alpha_{i,2} } \cdots S_{i,n_i} {O_0}^{\alpha_{i,n_i} } S_i^{*} x \\ & +
t \left(
 \sum_{k=1}^{n_i} (-1)^{\beta_{i,k} }\left( \prod_{j=1}^{k-1} S_{i,j} {O_0}^{\alpha_{i,j} } \right) S_{i,k} {O_0} ^{\beta_{i,k} }  \right) \mbox{,}
\end{split}
\end{equation}
where $\beta_{i,j} =0$ if $\alpha_{i,j}=1$ and $\beta_{i,j} =-1$ if $\alpha_{i,j}=-1$.

For every $k \in \{1 \stb n_i \}$ we define
\[ U_{i,k}=\left( \prod_{j=1}^{k-1} S_{i,j} {O_0}^{\alpha_{i,j}} \right) S_{i,k} \]
and let
 \[ \hat{U_i}=(-1)^{\beta_{i, n_i}} U_{i, n_i} {O_0}^{\beta_{i, n_i}} \mbox{ and } \mathring{U}_i = (-1)^{\beta_{i, 1}} S_{i,1} O_0^{\beta_{i,1}}\mbox{.}\]
 Using the previous notation one can see from equation \eqref{hosszu} that
\begin{equation}\label{u}
 U_i =S_{i,1} {O_0}^{\alpha_{i,1}} S_{i,2} {O_0}^{\alpha_{i,2} } \cdots S_{i,n_i} {O_0}^{\alpha_{i,n_i}} S_i^{*}
 \end{equation}
and we can also write
\[U_i= U_{i,n}\cdot {O_0}^{\alpha_{i,n_i}} S_i^{*} \mbox{.} \]
The vectors $b_1$ and $b_2$ can be written as $V_i t$, where
 \begin{equation}\label{b}
 V_{i} = \sum_{k=1}^{n_i}(-1)^{\beta_{i, k}} U_{i,k} {O_0} ^{\beta_{i, k}} \mbox{.}
 \end{equation}
 Equation \eqref{e2x} can be reformulated as follows
$$ (I-U_1)^{-1}V_1t= (I-U_2)^{-1}V_2t.$$
 By Lemma \ref{inverzaltalaban}, every entry of $(I-U_i)^{-1}$ is a rational function of the entries of $U_i$, which is generated by $O_0\stb O_m$.
Using Lemma \ref{ftlnremark} and the algebraic independence assumption on the coordinates of $t$ and $v_i$ it is clear that the previous equation holds for every vector $s \in \rr^d$ and $O'_1\stb O'_m\in SO(d)$. Thus we can eliminate $t$ from the previous equation and we get
\begin{equation}\label{etnelkul} (I-U_1)^{-1}V_1= (I-U_2)^{-1}V_2.
\end{equation}
First, we prove that it is enough to deal with the four dimensional case.

\begin{remark}\label{elokeszit}
\rm{Let us assume that $2s \ne 4$. From now on, we substitute block matrices into $O_i$ for $i=0,1 \stb m$, of the form
\[ M=\begin{pmatrix} N_1 & 0 \\ 0& N_2 \end{pmatrix} \mbox{,}\]
where $N_1 \in M(\theta) \subset SO(4)$ and $N_2 \in SO(2s-4)$. Since multiplying and adding these matrices we can count with the blocks separately.
 Clearly, a block matrix is invertible if and only if every block is invertible.

We need to guarantee that after the substitution, $I-U_1$ and $I-U_2$ are invertible. Since $2s-4 \ge 4$, Lemma \ref{nincsfixpont} shows that the group $SO(2s-4)$ contains a free subgroup (freely generated by $m$ elements) consisting of fix point free elements.
In order to prove that equation \eqref{etnelkul} does not hold for some substitution, it is enough to prove it for four dimensional matrices as in the following proposition.}
\end{remark}

\begin{proposition}\label{helyettesitosproposition}
We can substitute elements of the group $K$, defined in Lemma \ref{l40}, into $O_i$ for $i=0,1 \stb m$ such that equation \eqref{etnelkul} does not hold.
\end{proposition}
\proof

Substituting words of $A$ and $B$ we may assume that $U_1$ and $U_2$ are in $\mathcal{M}_1(\theta) \subset \mathcal{M}_1$.
Clearly, $U_i \in \mathcal{M}_1$ is of the following form for $i=1,2$:
$$
\left(
\begin{matrix}
p_i & -q_i & -r_i & -s_i\\
q_i & p_i & -s_i & r_i \\
r_i & s_i & p_i & -q_i \\
si & -r_i& q_i & p_i\\
\end{matrix}
\right).$$
Using Lemma \ref{inverz} and the fact that $W_1$ and $W_2$ are non-empty words, we get
\begin{equation*}
\begin{split}
(I-U_i)^{-1}&=\frac{1}{2-2 p_i}\left(
\begin{matrix}
1-p_i & q_i & r_i & s_i\\
-q_i & 1-p_i & s_i & -r_i \\
-r_i & -s_i & 1-p_i & q_i \\
-s_i & r_i& -q_i & 1-p_i\\
\end{matrix}\right) \\ &=\frac{1}{2-2p_i}(I-U_i)^T \end{split} \end{equation*}
since $U_i\in SO(4)$ and $(I-U_i)+(I-U_i)^{T}$ is a scalar matrix
and $p_1, p_2 \ne 1$. Equation \eqref{e2x} can be reformulated as
 $$\frac{1}{2-2p_1}(I-U_1)^T V_1= \frac{1}{2-2p_2}(I-U_2)^T V_2 \mbox{.}$$
This is equivalent to
\begin{equation}\label{e3}
(1-p_2)(I-U_1)^T V_1= (1-p_1)(I-U_2)^T V_2.
\end{equation}
Using equation \eqref{b} we get
\begin{equation}\label{e4}
\begin{split}
&(1-p_2)(I-U_1^T) (\rg{U}_1 + \sum_{k=2}^{n-1} (-1)^{\beta_{1,k} } U_{1,k} {O_0} ^{\beta_{1,k} } + \hat{U}_1) = \\
&(1-p_1)(I-U_2^T) ( \rg{U}_2 + \sum_{k=2}^{n-1} (-1)^{\beta_{2,k} } U_{2,k} {O_0} ^{\beta_{2,k} } +\hat{U}_2 ).
\end{split}
\end{equation}
Let
\begin{align}\label{align}
&M_1=(I-U_1^T) (\rg{U}_1 + \sum_{k=2}^{n-1} (-1)^{\beta_{1,k} } U_{1,k} {O_0} ^{\beta_{1,k} } + \hat{U}_1) \nonumber \\
\mbox{and similarly} \\
&M_2=(I-U_2^T) ( \rg{U}_2 + \sum_{k=2}^{n-1} (-1)^{\beta_{2,k} } U_{2,k} {O_0} ^{\beta_{2,k} } +\hat{U}_2 ) \mbox{.} \nonumber \end{align}
\begin{remark}\label{leghosszabb}
\rm{Equation \eqref{e4} depends only on the matrices $O_0, O_1 \stb O_m$. For an element of $O \in H=\langle O_0 \stb O_m\rangle$ we denote by  $\overline{O}$ the element of $K$ what we get after the substitution. Since $\ov{M}_1$ and $\ov{M}_2$ are generated by $A$ and $B$ defined in Lemma \ref{l40} we can write $\ov{M}_1=M_{\theta}(p_1,q_1,r_1,s_1)$ and $\ov{M}_2=M_{\theta}(p_2,q_2,r_2,s_2)$.
By expanding the brackets in equation \eqref{align} we get a sum where every summand is a subword or the inverse of a subword of $U_1$ and $U_2$ endowed with a sign. It is easy to see from Observation \ref{dot} \ref{dot4} that in order to determine the degree of the matrix in equation \eqref{align} we have to find the longest summands after the substitution. Basically the longest subword and the longest inverse of a subword occurring in $M_i$ are $\hat{U}_i$ and $-U_i^{T} \rg{U}_i =-U_i^{-1} \rg{U}_i$, respectively.}
\end{remark}
From now on we distinguish five major cases:
\begin{enumerate}[(a)]
\item\label{casea} $F_1 = F^{-1}$
\item\label{caseb} $F_1=F$, and there is no more $F^{\pm 1}$ on the paths $(P,Q)$.
\item\label{casec} $F_1=F$ and $F_2=F$
\item\label{cased} $F_1=F$ and $F_2=F^{-1}$.

Some of these cases originate in case \ref{casea}.

Case \ref{caseb} $\Rightarrow$ Case \ref{casea}: If there is only one $F^{\pm 1}$ on the path $(P,Q)$, then we just change the role of $P$ and $Q$ and we get case \ref{casea}.

Case \ref{casec} $\Rightarrow$ Case \ref{casea}: If $F_1=F$ and $F_2=F$, then we can change the role of the paths such that the common path of $W_1$ and $W_2$ contains $F_2$ instead of $F_1$.
This is again case \ref{casea}.

However, case \ref{cased} does not originate in case \ref{casea}, we can modify it to get a simpler form. In this case the role of $F_1$ and $F_2$ is symmetric hence we may assume that $F_1$ is not further from $P$ than $F_2$. This implies that there are some $O_i^{\pm 1}$'s on the path from the head of $F_2$ to the tail of $F_1$ (see figure Case \ref{casee}) which are not on the path $(P,Q)$ since the letters $F$ and $F^{-1}$ cannot succeed each other on a cycle.

 Thus, instead of to case \ref{cased} it is enough to investigate the following case:

\begin{center}
\includegraphics[width=6cm]{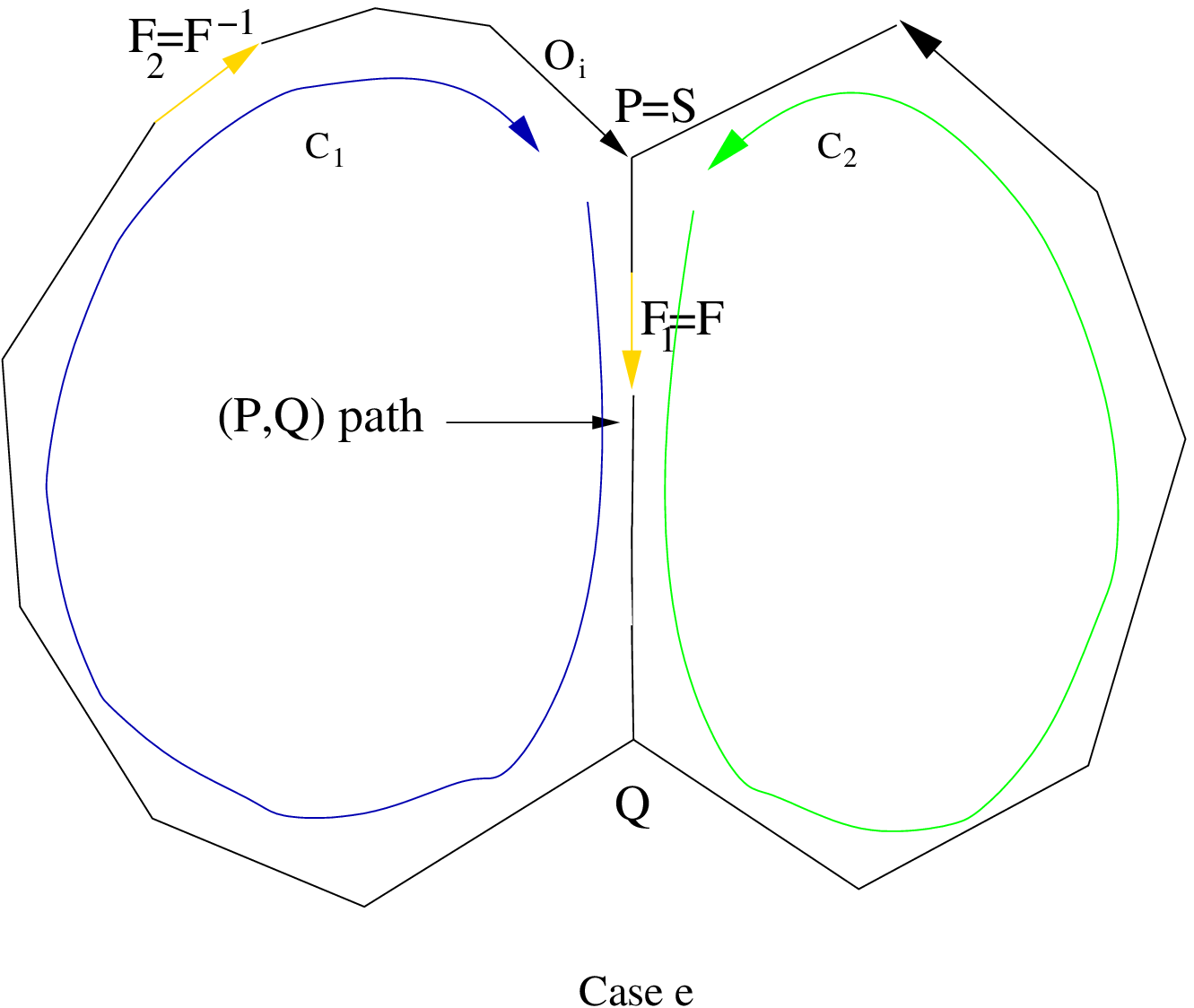}
\end{center}

 \item\label{casee} $F_1=F$ and $F_2=F^{-1}$ and there are some $O_i$'s in the path $(H(F_2),T(F_1))$ which are not in $(P,Q)$.
 Moreover we may assume that $F_1$ is not the only $F$ or $F^{-1}$ on its way since we assume that this case does not originate in case \ref{casea}. Using the same argument, the last $F^{\pm 1}$ on the $(P,Q)$ path has to be $F^{-1}$, otherwise we change the role of $P$ and $Q$. By symmetry again, we may assume that $F_1$ is the closest $F$ or $F^{-1}$ to $P$ on $\Delta$.
  \end{enumerate}

There are two major cases left and in both cases (and in every subcase) the starting point will be $S=P$.
Clearly, $S_1^{*} = S_2^{*}$ in this case so we denote it by $S_{*}$.

Now we substitute $O_i=A^{i_1 } B^{i_3} A^{i_2}$ for $i= 1\stb m$, where the absolute value of the exponents are pairwise different integers and $|i_3| \ge |i_1|,|i_2|>1$. Further, according to the case we investigate we substitute $O_0= \overline{S}_{*} A^{\varepsilon_1}\cdot B^D\cdot A^{\varepsilon_2} $ or $O_0= A^{\varepsilon_1}\cdot B^D\cdot A^{\varepsilon_2} \overline{S}_{*}^{-1} $, where $D, \varepsilon_1$ and $\ep_2$ will be chosen later.

We claim that if $D$ is large enough, then this substitution is monotone. More precisely we have the following.
\begin{lemma}\label{monotone}
\begin{enumerate}[(a)]
\item\label{monotonea}
Let $O_i=A^{i_1 } B^{i_3} A^{i_2}$ for $i= 0\stb m$, where the absolute value of the exponents are pairwise different integers and $|i_3| \ge |i_1|,|i_2|$.
Let $V$ be a reduced word on the letters $O_0 \stb O_m$ and let $V=U O_i^{\pm 1}$ or $V=O_i^{\pm 1} U$ and let $U$ be a subword of $V$. Then $\lg(\ov{U}) <\lg(\ov{V})$.
\item\label{monotoneb}
Let $O_i=A^{i_1 } B^{i_3} A^{i_2}$ for $i= 1\stb m$, where the absolute value of the exponents are pairwise different integers and $|i_3| \ge |i_1|,|i_2|>1$.
Let $\sigma_i^{1}$ and $\sigma_i^{2}$ denote the sum of the absolute value of the exponent of $O_i$'s occurring in $U_1$ and $U_2$, respectively. Let us assume that $|\ep_1| = |\ep_2|=1$ for every $1 \le i \le m$ and let
\[ D=\sum_{i=1}^{m} (\sigma_i^{1}+\sigma_i^{2})(|i_1|+|i_3|+ |i_2|) \mbox{.} \]
Then for every pair of subwords $V_1$ and $V_2$ of $W_j^{\pm 1}~(j=1,2)$ we have $\lg(\overline{V_1}) > \lg(\overline{V_2})$ if $V_1$ contains more $O_0^{\pm 1}$ than $V_2$.
\end{enumerate}
\end{lemma}
\proof
\begin{enumerate}
\item
 It is enough to prove it, when $V$ is of the form $UO_i^{\pm1 }$ or $O_i^{\pm1 }U$. Since the absolute value of the exponents are different $B$ cannot be eliminated after the substitution so it is easy to see that $\lg(\ov{U}) < \lg(\ov{U O_i^{\pm 1}})$ and $\lg(\ov{U}) < \lg(\ov{O_i^{\pm 1} U})$.
\item
It is easy to show that $B^{\pm D}$ cannot be eliminated from $\ov{O}_0^{\pm 1}$ since $|i_1|, |i_2|>1$ are different numbers.
\end{enumerate}
\qed

Case $\ref{casea}$:
\begin{center}
\includegraphics[width=5.8cm]{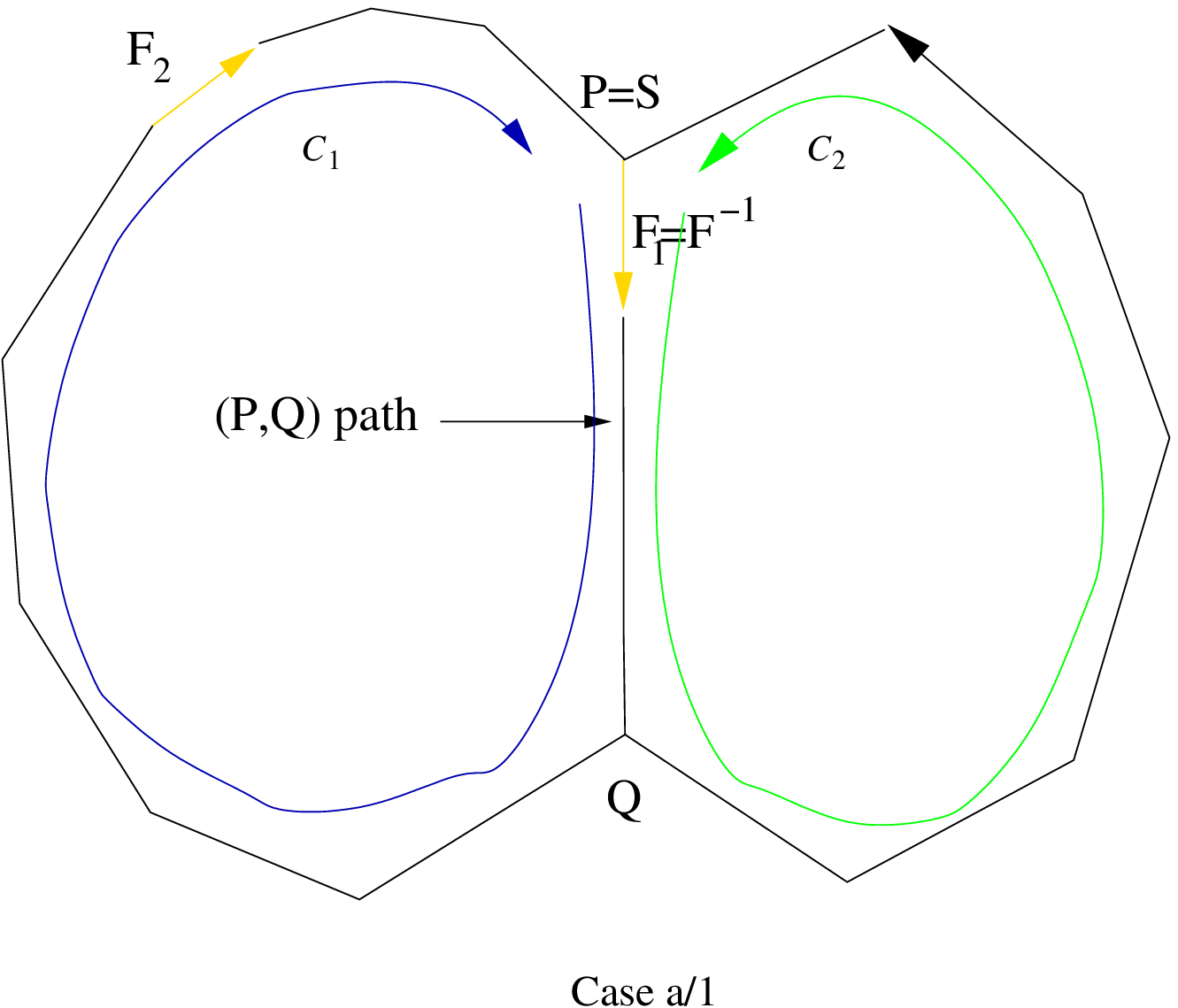}~
\includegraphics[width=5.8cm]{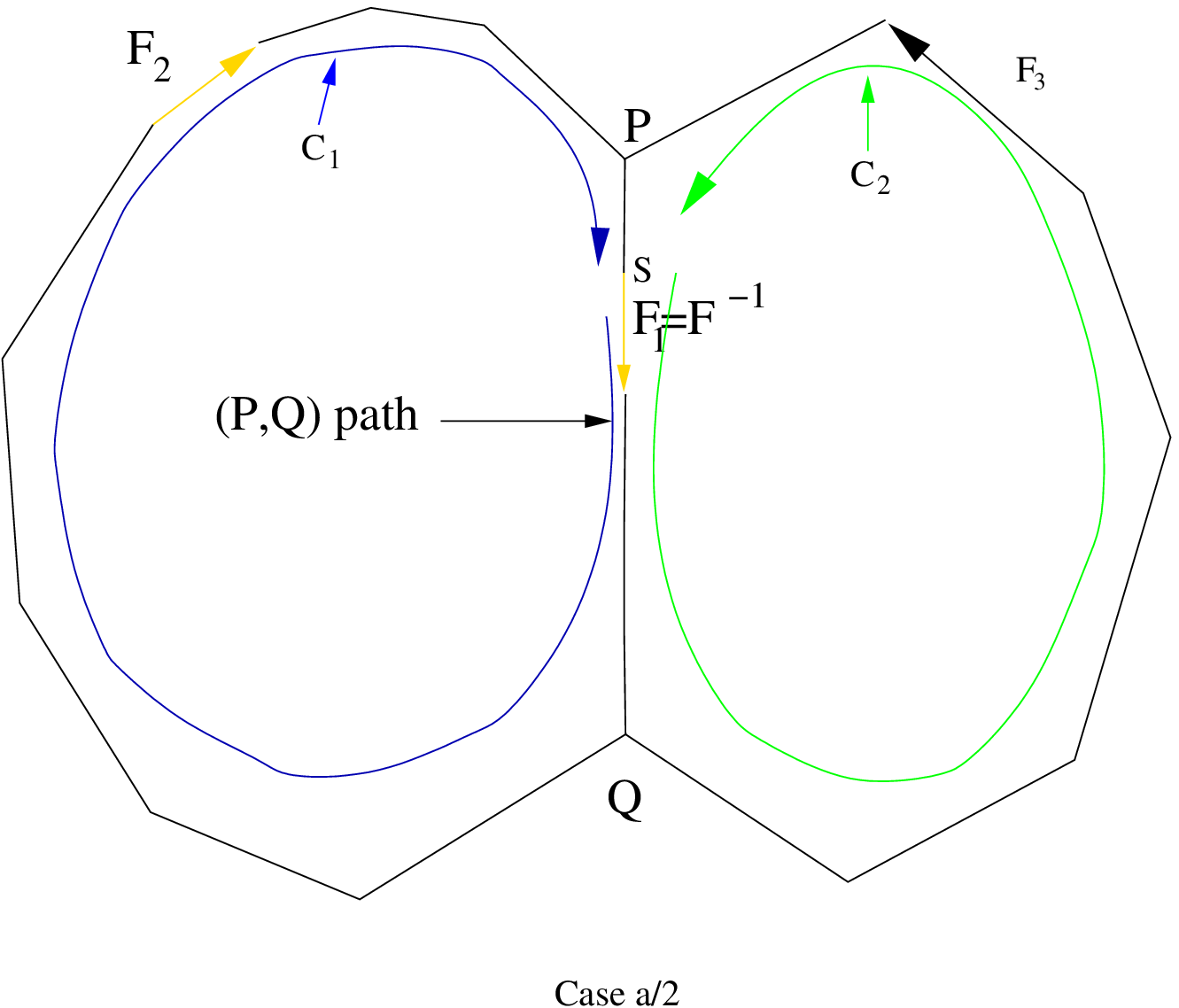}
\end{center}

\begin{enumerate}
\item
Let us assume first that the tail of $F_1$ equals to $P$.

Since the orthogonal transformations acts from the left, $W_1[1] \ne W_2[1]$ and neither of these letters are $F$ since both $W_1$ and $W_2$ represent a cycle.
Now we substitute $O_i=A^{i_1 } B^{i_3} A^{i_2}$ for $i= 0\stb m$, where the absolute value of all of these exponents are pairwise different integers and $|i_3| \ge |i_1|,|i_2|$.

One can see from Lemma \ref{monotone} \ref{monotonea} that the longest summands in $\ov{M}_1$ and $\ov{M}_2$ are $-\ov{U}_1$ and $-\ov{U}_2$, respectively since $\hat{U}_1 =-U_1$ and $\hat{U}_2 =-U_2$, while the first letter of $W_i$ is not $F$ so $\rg{U}_1, \rg{U}_2 \ne 1$ and every other summand in equation \eqref{e4} is a subword of one of them. By Observation \ref{dot} \ref{dot4} we have $\deg(\ov{U}_i) = \deg(\ov{M}_i)$ for $i=1,2$.

 Both $-\ov{U}_1$ and $-\ov{U}_2$ starts and ends with $A$ or $A^{-1}$. Since $W_1[-1] = W_2[-1]$ we have $\overline{U}_1 [-1] =\overline{U}_2 [-1]$ and since $W_1[1] \ne W_2[1]$ we can choose $i_1, i_2~(0 \le i \le m)$ such that $\overline{U}_1[1] \ne \overline{U}_2 [1]$.
Thus by the symmetry of $U_1$ and $U_2$ we may assume $\overline{U}_1[1] \ne \overline{U}_1 [-1]$ and $\overline{U}_2[1] = \overline{U}_2 [-1]$.
Lemma \ref{foatlo} shows that $\deg(p_1) < \deg(\ov{M}_1)$ and $\deg(p_2) = \deg(\ov{M}_2)$. This gives that $(1-p_2)\ov{M}_1  \doteq (1-p_1)\ov{M}_2$ does not hold so $(1-p_2)\ov{M}_1  \ne (1-p_1)\ov{M}_2$.

It is important to note that what we proved here is that both $\ov{M}_1$ and $\ov{M}_2$ have a unique longest summand and exactly one of these summands has a shorter conjugate. These facts guarantee that $(1-p_2)\ov{M}_1 \ne (1-p_1)\ov{M}_2$.
\item
Let us assume that the tail of $F_1$ is not $P$.

We use Lemma \ref{monotone} to calculate
\begin{equation}\label{emd}{\deg}\left( (I-U_i^T)(\hat{U}_i+ \sum_{k=1}^{n_i-1} (-1)^{\beta_{i,k}} U_{i,k} {O_0} ^{\beta_{i,k} }) \right) \mbox{.}
\end{equation}

Now we substitute $\ov{O}_0 = \ov{S}_{*} A^{\ep_1} B^D A^{\ep_2} $. Since $F_1=F^{-1}$, we have
\begin{equation*} \begin{split} \ov{\hat{U}}_i &= -\ov{S}_{i,1} \ov{O}_0^{\alpha_{i,1}} \ov{S}_{i,2} \ldots \ov{S}_{i,n_i}  O_0^{-1}  \\&= -\ov{S}_{i,1} (\ov{S}_{*} A^{\ep_1} B^D A^{\ep_2})^{\alpha_{i,1}} \ov{S}_{i,2} \ldots \ov{S}_{i,n_i}  A^{-\ep_2} B^{-D} A^{-\ep_1} {\ov{S}_{*}}^{-1} \mbox{.}
 \end{split} \end{equation*}

If $\alpha_{i,1}=-1$, then $\ov{U_i^{-1} \rg{U}}_i $ contains less $O_0^{\pm 1}$ than $\ov{\hat{U}}_i$.

If $\alpha_{1,1}=1$ (i.e. $F_2=F$), then
\begin{equation*} \begin{split}
\ov{{U}_i^{-1} \rg{U}_i} &= {\ov{S}_{*}}^{-1} \ov{S}_{*} A^{\ep_1}B^D A^{\ep_2} {\ov{S}_{i,n_i}}^{-1} \ldots \ov{S}_{i,2}^{-1}
 A^{-\ep_2} B^{-D} A^{-\ep_1} {\ov{S}_{*}}^{-1} \ov{S}_{i,1}^{-1} \ov{S}_{i,1} \\&= A^{\ep_1}B^D A^{\ep_2} {\ov{S}_{i,n_i}}^{-1} \ldots \ov{S}_{i,2}^{-1}
 A^{-\ep_2} B^{-D} A^{-\ep_1} {\ov{S}_{*}}^{-1}
 \end{split}
\end{equation*}
and
\[ \ov{\hat{U}}_i= -\ov{S}_{i,1} \ov{S}_{*} A^{\ep_1} B^D A^{\ep_2} \ov{S}_{i,2} \ldots \ov{S}_{i,n_i}  A^{-\ep_2} B^{-D} A^{-\ep_1} \ov{S}_{*}^{-1} \mbox{.}\]

It is easy to see that $S_{*} S_{i,1} \ne e$ since the path corresponding to $S_{*} S_{i,1}$ on $\Delta$ is non-trivial. It follows that $\ov{S}_{*} \ov{S}_{i,1}  \ne e$, which is equivalent to $ \ov{S}_{i,1} \ov{S}_{*} \ne e$. This shows using Lemma \ref{monotone} \ref{monotoneb} as well that $\ov{\hat{U}}_i$ is the longest summand of $\ov{M}_i$ again.

By the symmetry of paths between $P$ and $Q$ we may assume that if $\alpha_{i,1}=1$ for $i=1$ or $2$, then $\lg(S_{i,1}) \ge \lg(S_{*})$. This gives that $S_{i,1} S_{*}[1]= S_{i,1}[1] \ne e $ if $\alpha_{i,1}=1$.

We may assume that $(P,F^{-1}(P)) \notin E(\Delta)$ so if $S_{i,1} = e$ for $i=1 \mbox{ or }2$, then $\alpha_{i,1} =-1$ since otherwise this case goes back to case \ref{casea}. This also implies that $S_{1,1} = e$ and $S_{2,1} = e$ cannot hold at the same time.

Let us assume that neither $S_{1,1}$ nor $S_{2,1}$ is $e$. We also have
$S_{1,1}[1] \ne S_{2,1}[1]$ since the corresponding paths end in $P$. Therefore for suitable choice of the sign of the exponents $i_1$ and $i_2$ we may assume that $\ov{S}_{1,1}[1] = A^{e_1}=\ov{\hat{U}}_1[1]$ and $\ov{S}_{2,1}[1] = A^{e_2}=\ov{\hat{U}}_2[1]$
with $e_1 e_2<0$. It is easy to see that $\ov{\hat{U}}_1[-1] =\ov{\hat{U}}_2[-1]$ since $S_{*} \ne e$ so exactly one of $\ov{\hat{U}}_1$ and $\ov{\hat{U}}_2$ has shorter conjugate.

Let us assume that $S_{1,1}=e$. We have already proved that $S_{2,1} \ne e$ and $\alpha_{1,1}=-1$ in this case. Then $\ov{\hat{U}}_1[1]=A^{-\ep_2}$. Since $\hat{U}_1[-1], \hat{U}_2[1]$ and $\hat{U}_2[-1]$ are in $H'$ we have that for any choice of $i_1, i_2$ for $i=1 \stb m$ we may choose $\ep_2$ such that exactly one of $\ov{\hat{U}}_1$ and $\ov{\hat{U}}_2$ has shorter conjugate.

Similar result can be proved if $S_{2,1}=e$ so for suitable substitution we have $(1-p_2)\ov{M}_1  \doteq (1-p_1)\ov{M}_2$.
\end{enumerate}
Case \ref{casee}: Let us assume that this case does not originate in case \ref{casea}. It implies that $\alpha_{1,1}=\alpha_{2,1}=-1$ if $F_3$ is not on the $(P,Q)$ path. We have already assumed that $F_1$ is not the only $F^{\pm 1 }$ on the $(P,Q)$ path and the last $F^{\pm 1 }$ is $F^{-1}$ so if $F_3$ is on $(P,Q)$ path, then $F_3=F^{-1}$ again.

Therefore
\[ U_i=S_{i,1} {O_0}^{-1} S_{i,2} {O_0}^{\alpha_{i,2} } \cdots S_{i,n_i} O_O S_i^{*} \mbox{.}\]
As in the previous cases we write $S_{*}=S_1^{*} =S_2^{*}$.
Again we substitute $O_i=A^{i_1 } B^{i_3} A^{i_2}$ for $i= 1\stb m$, where the absolute value of the exponents are pairwise different integers and $|i_3| \ge |i_1|,|i_2|$ and $O_0=A^{\ep_1} B^{D} A^{\ep_2} \ov{S}_{*}^{-1}$, where $D$ is as large as in Lemma \ref{monotone} \ref{monotoneb}. By Lemma \ref{monotone} \ref{monotoneb} there are two possible choices for the longest term in $\ov{M}_1$. One of them is
\[ \ov{U}_1^{-1} \ov{S}_{1,1} \ov{O}_0^{-1} =\ov{S}_{*}^{-1}  \ov{O}_0^{-1}  \ov{S}_{1,n_1}^{-1} \ldots \ov{S}_{1,2}^{-1} \ov{O}_0 \ov{S}_{1,1}^{-1} \ov{S}_{1,1} \ov{O}_0^{-1} \mbox{,} \]
which equals to
\[  \ov{S}_{*}^{-1} \ov{S}_{*} A^{-\ep_2} B^{-D}  A^{-\ep_1} \ov{S}_{1,n_1}^{-1} \ldots \ov{S}_{1,2}^{-1} =A^{-\ep_2} B^{-D}  A^{-\ep_1} \ov{S}_{1,n_1}^{-1} \ldots \ov{S}_{1,2}^{-1} .
\]
The other one is
\[ \ov{S}_{1,1}  \ov{O}_0^{-1} \ov{S}_{1,2} \ldots \ov{S}_{1,n_1} =\ov{S}_{1,1} \ov{S}_{*} A^{-\ep_2} B^{-D} A^{-\ep_1} \ov{S}_{1,2} \ldots \ov{S}_{1,n_1} .\]

Since both $S_{1,1}^{-1}$ and $S_{*}$ starts at $P$ we have $S_{1,1}^{-1} \ne S_{*}$ hence \[ L_{1}=\ov{S}_{1,1} \ov{S}_{*} A^{-\ep_2} B^{-D} A^{-\ep_1} \ov{S}_{1,2} \ldots \ov{S}_{1,n_1}\] is the longest term of $\ov{M}_1$.
Similarly, the longest term of $\ov{M}_2$ is
\[ L_2 = \ov{S}_{2,1} \ov{S}_{*} A^{-\ep_2} B^{-D} A^{-\ep_1} \ov{S}_{2,2} \ldots \ov{S}_{2,n_2}.\]
We have already also assumed that $F_1$ is not the only $F$ or $F^{-1}$ on the $(P,Q)$-path so $\ov{S}_{2,n_2} =\ov{S}_{1,n_1}$.
Further, $S_{1,1}[1]$ and $S_{2,1}[1]$ are different since their tail is $P$ and.
The assumption that $F_1$ is the closest to $P$ among $F_1, F_2, F_3$ shows that $L_i[1]=\ov{S_{i,1} S_{*}}[1]=\ov{S}_{i,1}[1]$                                          for $i=1,2$.
Therefore we may choose the exponents $i_1$ and $i_2$ such that exactly one of $L_1$ and $L_2$ has shorter conjugate in $K$. This gives that for exactly one of $\ov{M}_1$ and $\ov{M}_2$ takes its degree in the main diagonal, which gives again $(1-p_2)(\ov{M}_1) \ne (1-p_1)(\ov{M}_2)$, finishing the proof of Proposition \ref{helyettesitosproposition}.

\qed

\section{Construction of the congruent pieces in high dimension}\label{s6}
\subsection{A set of symmetries of the regular simplex}

In this section we select isometries satisfying the conditions of Lemma \ref{l1}. Therefore, this set of isometries gives a decomposition of the balls (either open or closed) in $\rr^{d}$, where $d=4$ or $d=2s$ with $s\ge 4$.

We denote by $|v|$ the standard Euclidean norm of a vector $v \in \rr^{d}$ and we use the induced norm $||M||= \sup_{v \ne 0} \frac{|Mv|}{|v|}$ for $M \in SO(d)$.
\begin{remark}\label{dense}
Suppose that $\phi_i\in SO(d)$ $(i=1\stb k)$ are orthogonal transformations. Then for every $\varepsilon>0$ and for every $i \in \mathbb{N}$ there exists $O_{i,j} \in SO(d)$ for $j=1 \stb l$ such that $||\phi_i-O_{i,j} ||<\ep$ and the matrices $O_{i,j}$ are independent.
\end{remark}
\proof It is easy to see that the parametrization $\alpha_d: \Omega \rightarrow SO(d)$ is a continuous function of $\omega \in \Omega$, where the topology on $SO(d)$ is defined by the induced norm. There exists an everywhere dense subset of $\Omega$ whose elements are algebraically independent over $\mathbb{Q}$, finishing the proof of the Remark \ref{dense}.
\qed

\begin{defin}\label{affin}
Let $A_1 \stb A_{d+1}$ be the vertices of a regular simplex $S_d$ such $A_1 \stb A_{d+1}$ are in the boundary of the unit ball $B_d$. For $k=1, \ldots ,d+1$, let $H_k$ denote the affine hyperplane containing $A_i$ for every $~ i\ne k$. Let $\mathcal{A}'_k$ denote the set of $A_i$ which is contained by $H_k$. For instance, $\mathcal{A}'_1=\{A_2 \stb A_{d+1}\}$.
\end{defin}
Let $O$ denote the origin of the unit ball $B_d$ of dimension $d\ge 2$.
It is easy to see that $H_k$ is perpendicular to the vector $\vect{OA_k}$.

\begin{lemma}\label{orto}
Let $H_1$ and $H_k$ be two affine hyperplanes as in Definition \ref{affin}. Then there is a $\phi_k\in SO(d), ~k=1 \stb d+1 $ such that $\phi_k(H_1)= H_k$. Furthermore $\phi_k(\mathcal{A}'_1)=\mathcal{A}'_k$.
\end{lemma}
\proof
It is enough to show that the statement is true for $k=2$. It is easy to check that there is a reflection $r \in O(d)$ which fixes the points $A_3\stb A_{d+1}$ and maps $A_1$ to $A_2$. Clearly, $r$ is not in $SO(d)$. Therefore, we take another reflection $r'$ which fixes the points $A_1, A_2 \stb A_{d-1}$ and maps $A_{d}$ to $A_{d+1}$. Then the composition $\phi_2 =  r \circ r'$  $SO(d)$ and $\phi_k(\mathcal{A}'_1)=\mathcal{A}'_2$.
\qed

The image $\mathcal{I}_x$ of $x \in \rr^d$ was defined in Section \ref{4k}. For the multiset $\mathcal{I}_x$ we write $\mathcal{I}_x \subset H \subset \rr^d$ if and only if every element of $\mathcal{I}_x$ is in $H$.
\begin{lemma}\label{core}
Let $\phi_k \in SO(d) ~ (k=2\stb d+1)$ as in Lemma \ref{orto} and we fix $\phi_1 =id$. Let $T_b(x)=x+b$, where $b\in \rr^d$ with $t=|b|<\frac{2}{3d+4}$. Then every point $x \in B_d$ has a preimage $y=\phi^{-1}_j(x)$ for some $j=1\stb d+1$ such that for every $z \in B(y, t) \cap B_d$ the multiset $\mathcal{I}_z \subset B_d$.
\end{lemma}
\proof
We write a vector $u \in  \rr^{d}$ as $u=(u_1, u_2 \stb u_d)$.

We may assume that $A_1 =(0,0, \ldots,0,1)$, where $A_1$ is a vertex of the simplex given in Definition \ref{affin} and the vector $b$ and $\vect{OA_1}$ have the same direction.
Since $\phi_1, \phi_2 \stb \phi_{d+1}$ are orthogonal transformations, in order to verify for some $z \in B_d$ that  $\mathcal{I}_z=\{\phi_1(x), \phi_2(x) \stb \phi_{d+1}(x), T_b(x)\} \subset B_d$ it is enough to verify that $T_b(z)=z+b \in B_d$.
It is easy to see that if $z\in B_d$ with $z_d<-\frac{t}{2}$, then $z+b \in B_d$.

Every affine hyperplane $H_k$ divides $B_d$ into two parts. Let $F_k$ denote the part containing the simplex $S_d$ and $E_k$ denote the other one.

 We denote by $B_d^{1-2t}=\{x\in \rr^d: |x|<1-2t\}$, $A_k^{1-2t}$ and $H_k^{1-2t}$  the objects what we get from $A_k$ and $H_k$ by contracting $B_d$ with ratio $1-2t$ from the origin 0, respectively. The affine hyperplane $H_k^{1-2t}$ divides $B_d$ into two parts. We denote by $F^{1-2t}$ and $E^{1-2t}$ the two parts of $B_d$ which contains and which does not contain the contracted simplex, respectively.

If $x \in B_d^{1-2t}$, then we choose $\phi_1=id$ so $y=x$.
It is easy to see that if $z \in B(y,t)\cap B_d=B(x,t)\cap B_d$, then  $z+b \in B_d$.

Since the average of the coordinates of the points of the simplex $S_d$ is $0$ we have that the last coordinate of $A_2 \stb A_{d+1}$ is $-\frac{1}{d}$.
And similarly, since the last coordinate of $A_1^{1-2t}$ is $1-2t$, the last coordinate of $A_k^{1-2t}$ ($k=2 \stb d+1$) is $-\frac{(1-2t)}{d}$.

If $x \in E_k^{1-2t}$, where $k$ is not necessarily unique, then we choose $\phi_k$. Lemma \ref{affin} gives that $y \in E_1'$. Since the last coordinate of the point in $E_1^{1-2t}$ is smaller than $-\frac{(1-2t)}{d}$ and $t=|b|<\frac{2}{3d+4}$, we have $z_d \le -\frac{(1-2t)}{d} +t <-\frac{t}{2}$ for every $z \in B_d \cap B(y,t)$. Therefore, $z+b \in B_d$.

Clearly, $ B_d \subseteq \bigcup_{k=1}^{d+1}  E_k^{1-2t} \cup B_d^{1-2t} $, finishing the proof of Lemma \ref{core}.
\qed
\begin{remark}\label{kulonbozokep}
\rm{For $x \in B_d$, in order to find a preimage $y$ which is a core point we use Lemma \ref{core} and besides, we guarantee that the elements of
$\mathcal{I}_y$ are different.}
\end{remark}
\subsection{Proof of Theorem \ref{t12} }\label{ss2}
Now we can complete the proof of Theorem \ref{t12}.

We construct $4(2s+1)+1$ maps satisfying the conditions
of Lemma \ref{l1}, where the bijections are orthogonal transformations of $SO(d)$ with $d=2s\ge 4$ and $d\ne 6$. In this case the $d$ dimensional unit ball
can be decomposed into $4(2s+1)+2$ pieces.

According to Lemma \ref{core}, there exist orthogonal transformations $\phi_i\in SO(d)$ ($i=1\stb d+1$) and $|b|=t < \frac{2}{3d+4}$
which have the property that every point has a preimage $y$ such that for every $z \in B(y,t)$ we have $\mathcal{I}_z \subset B_{d}$.

By Lemma \ref{dense}, for every $1 \le i \le 2s+1$ there exist independent orthogonal transformations $O_{i,1}, O_{i,2}, O_{i,3}$ and $O_{i,4}$
such that $||O_{i,j}-\phi_i||<t$ and $||O_0-id||<t$. Furthermore, we assume that $O_{i,j}$, $O_0$  and $b$ form an independent system with respect to the standard basis of $\rr^d$.

By choosing a suitable orthonormal basis in $\rr^d$, we may use Lemma \ref{core}. If $x \ne x_0 =0$, then
there exists an $i$ such that for every $j=1 \stb 4$ we have $\mathcal{I}_y \subset B_{d}$ for $y=O_{i,j}^{-1}(x)$. By Remark \ref{kulonbozokep}, if we can guarantee that $\mathcal{I}_y$ consists of different points, then $y$ is a core point.

Lemma \ref{nincsfixpont} shows that $O_{i_1,j_1}(y)\ne O_{i_2,j_2}(y)$ if $y\ne 0$, which is the case since $x \ne 0$.
The only case which remains is that
\begin{equation}\label{nemcore}
O_{k,l}(y)=O_0(y)+b
\end{equation}
for some $y$ and $O_{k,l}$, where $ k\in \{1\stb d+1\}, l\in \{1\stb 4\}$.

\begin{lemma}\label{1os}
Let $O_{m,n} ~(m=1 \stb 2s+1\mbox{, } n=1 \stb 4)$, $O_0$ and $b$ an independent system in $\rr^d$, where $d=2s$. For every $x\ne 0$, there exists at most one pair of linear transformations $O_{i,j}$
and $O_{k,l}$ such that for the point $y=O_{i,j}^{-1}(x)$ the equation $O_{k,l}(y)=O_0(y)+b$ might be satisfied.
\end{lemma}
\proof
Let us assume that for some $0 \ne x \in B_d$ we have
 \[ O_{k_1,l_1}(y_1)=O_0(y_1)+b \quad \mbox{and} \quad O_{k_2,l_2}(y_2)=O_0(y_2)+b \mbox{,}\]
where $y_1 =O_{i_1,j_1}^{-1}(x)$ and $y_2 =O_{i_2,j_2}^{-1}(x)$.
Thus, we get
\begin{equation}\label{coreseg}
(O_{k_1,l_1}-O_0)(O_{i_1,j_1}^{-1}(x))=b \mbox{ and }
(O_{k_2,l_2}-O_0)(O_{i_2,j_2}^{-1}(x))=b \mbox{.}
\end{equation}
Using again Lemma \ref{nincsfixpont} we get that $O_0^{-1}O_{i,k}-I$ is invertible.
Therefore equation \eqref{coreseg} can be written in the form
$$x=O_{i_1, j_1}(O_0^{-1}O_{k_1,l_1}-I)^{-1}O_0^{-1}b \mbox{ and } x=O_{i_2, j_2}(O_0^{-1}O_{k_2,l_2}-I)^{-1}O_0^{-1}b.$$
Hence
$$O_{i_1, j_2}(O_0^{-1}O_{k_1,l_1}-I)^{-1}O_0^{-1}b=O_{i_2, j_2}(O_0^{-1}O_{k_2,l_2}-I)^{-1}O_0^{-1}b.$$
Since $O_{m,n}, O_0$ and $b$ form an independent system, we can eliminate $b$ from the previous equation, and we get the following:
\begin{equation}\label{rovidhelyettesites}
O_{i_1, j_1}(O_0^{-1}O_{k_1,l_1}-I)^{-1}O_0^{-1}=O_{i_2, j_2}(O_0^{-1}O_{k_2,l_2}-I)^{-1}O_0^{-1}\mbox{.}
\end{equation}
Using Lemma \ref{ftln}, we may substituting $O_0=id$ and we get
\begin{equation}\label{retek}
O_{i_1, j_1}(O_{k_1,l_1}-I)^{-1}=O_{i_2, j_2}(O_{k_2,l_2}-I)^{-1}\mbox{.}
\end{equation}
If $O_{i_1, j_1} =O_{i_2,j_2}$ or $O_{k_1, l_1} =O_{k_2,l_2}$, then it is clear from equation \eqref{retek} that $O_{i_1, j_1} =O_{i_2,j_2}$ and $O_{k_1, l_1} =O_{k_2,l_2}$.

Thus we can assume that $O_{i_1, j_1} \ne O_{i_2,j_2}$ and $O_{k_1, l_1} \ne O_{k_2,l_2}$.
Then we substitute such that $\ov{O}_{i_1, j_1} = \ov{O}_{i_2,j_2}$. This implies $\ov{O}_{k_1, l_1} = \ov{O}_{k_2,l_2}$.  Hence either $O_{i_1,j_1} = O_{k_1,l_1}$ and $O_{i_2,j_2} = O_{k_2,l_2}$, or $O_{i_1,j_1}=O_{k_2,l_2}$ and $O_{i_2,j_2} = O_{k_1,l_1}$.
\begin{enumerate}
\item
Let us assume first that $O_{i_1,j_1} = O_{k_1,l_1}$ and $O_{i_2,j_2} = O_{k_2,l_2}$. We shortly denote $O_{i_1,j_1}$ by $U$ and we substitute $O_{i_2,j_2}=U^2$. From equation \eqref{retek} we get
\[ U(U-I)^{-1} =U^2 (U^2-I)^{-1}.\]
This gives
\[ U^2-I =(U-I)(U+I)=(U-I)U, \]
which is a contradiction since $U-I$ is invertible by Lemma \ref{nincsfixpont}.
\item
Let us assume that $O_{i_1,j_1}=O_{k_2,l_2}$ and $O_{i_2,j_2} = O_{k_1,l_1}$.
Then we denote $U= O_{i_1,j_1}$ and we substitute $O_{i_1,j_1}=U^2$ again.
Similar calculation gives
\[ I=(U+I)U .\]
This gives $U^2+U-I$, which is a polynomial expression, contradicting Lemma \ref{ftln}.
\end{enumerate} This shows that
equation \eqref{rovidhelyettesites} holds if and only if $\{i_1, j_1\} = \{i_2, j_2\}$ and  $\{k_1, l_1\} = \{k_2, l_2\}$, finishing the proof of Lemma \ref{1os}.
\qed

\bigskip
For every $x\ne 0$ we have already found $O_{i,1}$, $O_{i,2}$, $ O_{i,3}$ and $O_{i,4}$ such that $\mathcal{I}_{y_j} \subset B_d$, where $y_j={O_{i,j}^{-1}(x)}$ for every $j=1 \stb 4$. By Lemma \ref{1os} at least three of $y_j$ is a core point satisfying Lemma \ref{l1} \ref{fofeltetel2}.

If $x_0$ is the origin it has to satisfy condition Lemma \ref{l1} \ref{fofeltetel1}. Therefore we need to guarantee that $F^{-1}(x_0)=F^{-1}(0)$  is a core point.
Indeed, if $y=F^{-1}(0)\ne 0$, then $O_{i_1,k_1}(y) \ne O_{i_2,k_2}(y)$ holds again by Lemma \ref{nincsfixpont}. Clearly, $F(y)=0$, thus $O_{k,l}(y)=F(y)$ cannot hold for any $k$ and $l$. Due to the choice of $b$ we have $\mathcal{I}_y \subset B_{d}$.

Since the matrices $O_0$, $O_{i,j}$ and $b$ is an independent system and the dimension $d=2s \ge 4$ and $s\ne 3$, we can use Theorem \ref{foprop}. Thus the graph $\Ga_{\mathcal{F} }$
has the property \eqref{coin} for $\mathcal{F} = \left\{ F, O_{m,n} \right\} ~ (m=1\stb 2s+1, n=1\stb 4)$, where $F=T_bO_0$.

We conclude that $\mathcal{F}$ satisfies the conditions of Lemma \ref{l1}, finishing the proof of Theorem \ref{t12}.

\qed

\begin{remark}
\begin{enumerate}
{\rm \item
The proof above gives a construction for $m=4\cdot(2s+1)+1=4d+5$ pieces in dimension $d=2s$.
We can easily obtain a construction for $m> 4d+5$, since we can add any finite number of orthogonal transformations with algebraically independent parameters to the already defined ones, which satisfy the conditions of Lemma \ref{l1}.
    \item
Most probably, this bound $4d+5$ is practically not the best but this construction of Section \ref{s6} cannot be modified without difficulties.}
\end{enumerate}
\end{remark}

\section{Decomposition in higher dimension}\label{s7}
In this section we prove Theorem \ref{t2}.\\
In \cite{KL} the authors proves the following theorem
\begin{theorem}The $3s$ dimensional ball can be decomposed into finitely many pieces for every $s \in \mathbb{Z}^{+}$.
\end{theorem}
This shows that there is a decomposition for $d=6$ and $d=9$. In order to prove Theorem \ref{t2}, by Theorem \ref{t12} it is enough to prove it when $d\ge 7$ is odd and $d \ne 9$. Such an integer can be written in the form $d=d'+3$. Then we write the elements $x$ of $\rr^d$ as $x=(y, z)$, where $y \in \rr^{d'}$ and $z \in \rr^3$, where $d'\ge 4$ is even and $d' \ne 6$ .

We shall recall some of the results of \cite{KL} for the 3 dimensional case using our notation.

The following lemma is essentially the same as \cite[Lemma 3.5]{KL}.
\begin{lemma}\label{l5}
Suppose that $O'_0, O'_1 \stb O'_m \in SO(3)$ and $b \in \rr^{3}$ form an independent system. Let $F=T_b O'_0$ and $\mathcal{F}=\{ F, O'_1 \stb O'_m \}$. If $\mathcal{C}$ is a cycle in $\Ga_{\mathcal{F}}$, then the corresponding word does not contain the letter $F$ or $F^{-1}$.
\end{lemma}

We remind that a cycle has distinct points aside from the first and the last vertices of it which coincide.
We refer to \cite[Lemma 4.1]{KL} which states the following.

\begin{lemma}\label{l4}
Suppose that $O_0', O_1' \stb O_m' \in SO(3)$ and $b\in \rr^3$ form an independent system. Let $\mathcal{F}=\{T_b O_0', O_1' \stb O_m'\}$.
Then $\Gamma_{\mathcal{F}}$ has property \eqref{coin}.
\end{lemma}
Finally, a version of Lemma 4.2. in \cite{KL} states following
\begin{lemma}\label{l6}
Suppose that $ O_1' \stb O_m' \in SO(3)$ are independent orthogonal transformations. Then for every $0 \ne x \in \rr^{3}$ there are at most two elements of the form $O_i'$ such that $I_{O_{i}'^{-1}(x)}$ does not consist of different points.
\end{lemma}

Our aim is to construct $d$ dimensional special orthogonal transformations satisfying the conditions of Lemma \ref{l1}.

Let $b_1$ be a vector in $\rr^{d'}$.
Let $\phi_1, \phi_2 \stb \phi_{d'+1} \in SO(d')$ as in Lemma \ref{core} with the additional assumption that $\vect{OA_1}$ and $b_1$ have the same direction. For every $1 \le i \le d'+1$ we choose $20$ orthogonal transformations $O_{i,j} ~ (j=1\stb 20)$ such that $|| O_{i,j}-\phi_i || \le \ep'$ for some $\ep'>0$. Let $O_0 \in SO(d')$ satisfies $||O_0- I_{d'}|| \le \ep'$, where $I_{n}$ denotes the $n$ dimensional identity matrix and let $F=T_{b_1} O_0$. We assume that $O_0, O_{i,j}$ and $b_1$ form an independent system.

Let $\phi_1', \phi_2', \phi_3', \phi_4' \in SO(3)$ as in Lemma \ref{core} and $b_2 \in \rr^3$. We assume again that for one of the points $A_1'$  of the $3$ dimensional simplex, the vector $\vect{OA_1'}$ and $b_2$ have the same direction. We denote by $1 \le j' \le 4$ the integer such that $j \equiv j' \mbox{ mod } 4$.
For every $1 \le j \le 20$ we choose $O_{i,j}'  \in SO(3)$ for $( 1 \le i \le d'+1)$ such that $||O_{i,j}' -\phi_{j'}'|| \le \ep'$ and let $||O_0'-I_3|| < \ep'$ and let $F'=T_{b_2}O_0'$. We assume again that $O_0', O_{i,j}'$ and $b_2$ form an independent system, where $b_2 \in \rr^3$.

 We define the orthogonal transformations $\hat{O_0}$ and $\hat{O}_{i,j} \in SO(d)$ by \[ \hat{O}_0(x,y)=\left(O_0(y),O_0'(z)\right) \mbox{ and } \hat{O}_{i,j}(y,z) = (O_{i,j}(y), O_{i,j}'(z)  ) \]
 and for $b=(b_1,b_2) \in \rr^d$ let
 \[ \hat{F} =  T_{b}\hat{O}_0(x) =(T_{b_1}O_0(y) , T_{b_2}O_0'(z) ) \mbox{.} \]
One can see that there exists some $\ep$ depending only on $\ep'$ and $d$ such that for $\Phi_{i,j}=(\phi_i,\phi_j') \in SO(d)~(i=1 \stb d+1, ~ j=1\stb 4)$ there are at least $5$ orthogonal transformations of the form $\hat{O}_{i,j}$ such that $|| \hat{O}_{i,j}- \Phi_{i,j}|| <\ep$ and if $\ep'$ tends to $0$ then $\ep$ tends to $0$ as well.

\begin{lemma}\label{utolsolemma}
 The graph $\Gamma_{\mathcal{\hat{F}}}$ has property $(\ref{coin})$, where $\mathcal{\hat{F}}= \left\{ \hat{F}, \hat{O}_{1,1}  \stb \hat{O}_{d'+1,20}   \right\}$.
\end{lemma}
\proof
We claim that if $\mathcal{C}$ is a cycle in $\Gamma_{\hat{F}}$, then for every vertex $(y,z) \in \rr^d$ on the cycle $\mathcal{C}$ we have $y=0$.
One can easily assign to $\mathcal{C}$ a word $W$ by identifying the vertices by the letters $\hat{O}_{i,j}$ and $\hat{F}$. Clearly, $W$ is a reduced word.

Let $W'$ be the restriction of $W$ to the last three coordinates. Lemma \ref{l5} shows that if $W'$ contains $F^{\pm 1}$, then $W'$ does not have a fixed point. Hence $W$ does not contain the letter $\hat{F}$ or $\hat{F}^{-1}$. In this case $W''$ which is the restriction of $W$ to the first $d'$ coordinates can be identified by an element of the group $K$ which is not the identity element. Lemma \ref{nincsfixpont} shows that the only fixed point of $W''$ is $0$ so $y=0$ for each vertex $(y,z)$ of $\mathcal{C}$.

Thus if $\Gamma_{\hat{F}}$ contains two cycles $\mathcal{C}_1$ and $\mathcal{C}_2$ sharing an edge, then $\mathcal{C}_1$ and $\mathcal{C}_2$ can be considered as cycles in $\rr^{3}$ in the graph $\Gamma_{\mathcal{G}}$, where $\mathcal{G}=\left\{ F', O_1' \stb O_m' \right\}$. In this case \cite[Lemma 2.1]{KL} shows that $\mathcal{C}_1$ and $\mathcal{C}_2$ coincide, finishing the proof of Lemma \ref{utolsolemma}.
\qed

\bigskip
We prove that $x_0 =0$ satisfies the conditions given in Lemma \ref{l1} with respect to the graph $\Gamma_{\mathcal{F}'}$.

It is easy to see that for $\hat{F}^{-1} (x_0)=\hat{F}^{-1} (0)=(y_0,z_0)$ we have $0 \ne y_0 \in \rr^{d'}$ which shows that $\hat{O}_{i,j}((y_0,z_0)) \ne \hat{O}_{k,l}((y_0,z_0))$ if $(i,j) \ne (k,l)$ using Lemma \ref{nincsfixpont}. Moreover these points differ from $\hat{F}(\hat{F}^{-1})(0) =0$ hence $\hat{F}^{-1}(0)$ is a core point if $|b| <1$.

\begin{lemma}\label{l7}
For every $0 \ne x=(y,z) \in \rr^{d}$ there are at most two $\hat{O}_{i,j}$ such that $\mathcal{I}_{\hat{O}_{i,j}^{-1}(x)}$ does not consist of different elements.
\end{lemma}
\proof
If $y \ne 0$, then Lemma \ref{1os} shows that there is at most one $\hat{O}_{i,j}$ such that the elements of $\mathcal{I}_{\hat{O}_{i,j}^{-1}(x)}$ are not different.
Lemma \ref{l5} shows that $O_{i,j}'(v) \ne T_{b_2}O_0'(v)$ for every $v \in \rr^{3}$ so it can only happen that for $(0,v)=\hat{O}_{i,j}^{-1}(x) \in \rr^d$ we have $\hat{O}_{k_1,l_1}(0,v)=\hat{O}_{k_2,l_2}(0,v)$ for some $1 \le k_1, k_2 \le d' $ and $1 \le  l_1,l_2\le 20$. Lemma \ref{l6} gives that there at most two $\hat{O}_{k,l}$ with this property.
\qed

One can see from Lemma \ref{l7} that in order to verify that every $0 \ne x=(y,z) \in B_d$ is the image of at least three core points we only have to show that for at least five $(i,j)$ pair $\mathcal{I}_{\hat{O}_{i,j}^{-1}(x)} \subset B_d$.

We will compare the length of vectors of different dimension so we denote by $|v|_k$ the length of the vector $v \in \rr^k$.
Now we assume that $|b_1|_{d'}=|b_2|_3 =t$.
Choosing a suitable basis again, we assume that $b=(b_1,b_2)$, where $b_1 =(0 \stb 0,0,t)\in \rr^{d'}$ and $b_2 =(0,0,t) \in \rr^{3}$ and $A_1 \in \rr^{d'}$, which is a vertex of the simplex defined in Lemma \ref{affin}, is just $\underbrace{(0 \stb 0,1)}_{d'}$.
\begin{lemma}\label{geom}
\begin{enumerate}[(a)]
\item\label{geom1}
Let $x \in \rr^k$ with $|x|_k>\frac{1}{3}$. Then there exists $\phi_i \in SO(k)$ as in Lemma \ref{affin} such that for $\phi_i(x)=(y_1 \stb y_k)$ we have $y_k \le -\frac{1}{3k}$.
Moreover, if $|c|_k=|(0 \stb 0,c_k)|_k \le \frac{1}{6k}$, then for every $u \in B(\phi_i(x),|c|_k)$ with $|u|_k=|x|_k$ we have $u_k \le -\frac{1}{6k}$ and $|u+c|_k \le |x|_k$.
\item\label{geom2}
For every $\ep>0$ there exists $r>0$ such that if $|c|<r$, then for every $x \in B_k \setminus B_k^{\ep}$ there exists $\phi_i \in SO(k)$ such that $|u+c|_k  \le |x|_k$ for every $u \in B(\phi_i(x),|c|_k)$ with $|u|_k =|x|_k=|\phi_i(x)|_k$.
\end{enumerate}
\end{lemma}
\proof
\begin{enumerate}[(a)]
\item
Clearly, $B_k \setminus B_k^{\frac{1}{3}} \subset \cup_{i=1}^{k+1} E_i^{\frac{1}{3}}$, where $E_i^{\frac{1}{3}}$ denotes the intersection of a half-plane with $B_k$ as in Lemma \ref{core}. Hence $x$ is contained in $E_i^{\frac{1}{3}}$ for some $1 \le i \le k+1$. Since $\phi_i(E_i^{\frac{1}{3}}) =E_1^{\frac{1}{3}}$, the last coordinate of $\phi_i(x)$ is smaller than or equal to $-\frac{1}{3k}$. Finally, one can easily verify that the last coordinate of $z$ is smaller than $-\frac{|c|_k}{2}$, which gives that $|u+c|_k \le |u|_k =|x|_k$.
\item
Using the same argument again, we may assume that $\phi_i(x) \in E_1^{\ep}$. This shows that the last coordinate of $u=(u_1 \stb u_k)$ is smaller than or equal to $- \ep \frac{1}{k}+|c|_k$. If $|c|_k$ is small enough, then $u_k < -\frac{|c|_k}{2}$ which guarantees that $|u+c|_k \le |u|_k$.
\qed
\end{enumerate}
\begin{lemma}\label{geomkieg}
If $\sqrt[d]{2t^2}=|b|_d=|(0 \stb 0,t,0,0,t)|$ is small enough, then for every $x=(y,z) \in \rr^{d}$ there exists $\Psi=(\phi_i, \phi_j')$ such that $u+b \in B_d$ for every $u =(u_1,u_2)\in B(\Psi(x), |b|)$ with $|u_1|_{d'}= |y|_{d'}$ and $|u_2|_3=|z|_3$.
\end{lemma}
\proof
It is easy to see that if $|x|_d < 1-|b|_d$, then $u+b \in B_d$.
Therefore, if $|b|_d <  \frac{1}{3}$, then we may assume that either $|y|_{d'}$ or $|z|_3$ is greater than $\frac{1}{3}$ since $\frac{\sqrt{2}}{3} + \frac{1}{3}<1$.

Let us assume first that $|y|_{d'}\ge \frac{1}{3}$.
Let $\theta_{d'}$ be a negative number what we will define later and let $\ep=-\frac{\theta_{d'}}{6}$. By Lemma \ref{geom} \ref{geom2} there exists $r>0$ such that if $|b_2|_3 <r$, then for every fixed $w =(w_1,w_2)\in \rr^{d'+3}$ with $|w_2|_3 > \ep$ there exists $\phi_j'$ such that for every $u_2 \in B(\phi_j'(w_2), |b_2|_3)$ with $|u_2|_3=|w_2|_3$ we have $|\phi_j'(u_2)+b_2|_3 \le |w_2|_3$. Therefore, if $|y|_{d'} \ge \frac{1}{3}$ and  $|z|_3 \ge \ep$, then there exists a $d$ dimensional orthogonal transformation of the form $\Psi=(\phi_i, \phi'_j)$ such that for every if $u =(u_1,u_2) \in \rr^d$ with $|u_1|_{d'}= |y|_{d'}$, $|u_2|_3= |z|_3$ and $u \in B(\Psi(x),|b|)$, then $u+b \in B_d$. Thus we may assume $|z|_3 < \ep$.

We show that $|u_1+b_1|_{d'}^2+|u_2+b_2|_3^2 \le |y|_{d'}^2+ |z|_3^2$.
This is equivalent to $|b_1|_{d'}^2+|b_2|_3^2 \le -2u_1b_1 -2u_2b_2$, where the product of two vectors is the standard inner product. Using $|b_1|=|b_2|=t$ we get
\begin{equation}\label{tau} t \le -|u_1|_{d'} \cos\tau_1- |u_2|_3 \cos\tau_2 \mbox{,} \end{equation}
where $\tau_1$ and $\tau_2$ denotes the angle between $b_1$ and $u_1$ and between $b_2$ and $u_2$, respectively.
Lemma \ref{geom} \ref{geom1} gives that the last coordinate $u_1$ is smaller than $-\frac{1}{6d'}$ so $\cos\tau_1$ can be estimated from above by a number $\theta_{d'}=-\frac{1}{3} \frac{1}{6d'}$ which only depends on $d$.
Thus \begin{equation}\label{teta} -|u_1|_{d'} \cos \tau_1- |u_2|_3 \cos \tau_2 \ge -\frac{1}{3} \theta_{d'}- |u_2|_3 \ge -\frac{1}{3} \theta_{d'} + \frac{1}{6} \theta_{d'} =-\frac{1}{6} \theta_d \mbox{.} \end{equation}
It is easy to see that this last term in equation \eqref{teta} is a positive number which only depends on $d$ so for suitable choice of $t$ combining with the previous conditions for $|b| = \sqrt[d]{2t^2}$ we have that equation \eqref{tau} holds, finishing the proof of Lemma \ref{geomkieg}.
\qed

Lemma \ref{utolsolemma}, Lemma \ref{l7} and Lemma \ref{geomkieg} imply that for every $x \in B_d$ there are at least five $\ov{O}_{i,j}$ such that $\mathcal{I}_{\ov{O}_{i,j}}^{-1}(x) \subset B_d $
and clearly, three of them are core points. We conclude that all the conditions of Lemma \ref{l1} are satisfied, finishing the proof of Theorem \ref{t2}.
\qed

Now we collect the results on the number of pieces required for the decomposition in different dimensions to prove Theorem \ref{t3}.

{\bf Proof of Theorem \ref{t3}}: \\
We distinguish 3 major cases.
\begin{enumerate}

\item If $d=3,6$ or $9$ then $B_d$ can be decomposed into finitely pieces by \cite[Theorem 1]{KL}.

\item If $d=2s \ge 4$, where $s\ge 2$ and $s \ne 3$, then the number of orthogonal transformations is $4(2s+1)+1=4d+5$ and hence $B_d$ can be decomposed into $(4 d+5)+1=4d+6$ pieces, by Theorem \ref{t11}.

\item If $d=2s+3$, where $s\ge 2$ and $s \ne 3$, then the number of orthogonal transformations is $5(4(2s+1)) +1=20(d-2)+1$ and hence $B_d$ can be decomposed into $(20(d-2)+1)+1=20d-38$ pieces, by Theorem \ref{t2}.
\end{enumerate}
Therefore the number of pieces is asymptotically, less than or equal to $20d$ if d.

\section{Problems and results in dimension $d$}\label{vegsos8}

In dimension $d=2$ the transformation group $O_2$ does not contain noncommutative
free subgroups, thus the methods worked  out in \cite{KL} and in this article and  cannot say anything about the divisibility of the
discs. C. Richter posed a question about decomposition of the disc using affine transformations instead of orthogonal transformations.
A celebrated result of von Neumann shows that the group of affine transformations contains noncommutative
free subgroups. In this case, the main difficulty is to satisfy the conditions \ref{fofeltetel2} of Lemma \ref{l1}. We do not know whether or not Richter's problem can be solved along these lines.

By Theorem \ref{t11}, the minimal number of pieces $\tau_d$ which is needed to decompose $B_d$, is less than $20 d$ for $d \ge 10$.
The main result of \cite{W1}, which was reproved in \cite{vD}, shows that $\tau_d>d$. Thus we get that $\tau_d=\Theta(d)$, which is best possible in some sense.
This widely improves the upper bound of $\tau_d$ for $d=3s$ given in \cite{KL}, where it was shown that $\tau_d \le \exp(c_1 d \log d)$ for a positive constant $c_1$.

As for $d=3,$ the question whether or not $B_3$ is $m$-divisible for
$4\le m\le 21$ is open, for $d\ge 10$, the question whether or not $B_d$ is $m$-divisible for $d+1\le m \le 20d$ also remains open.
There are several obstacles in the way
of improving these bounds. One of them is the condition of
Lemma \ref{l1} which requires that every point
$x\ne x_0$ has to be the image of at least three core points.
In \cite[Example 6.1]{KL} was shown that this condition of Lemma \ref{l1} is sharp.

The most related question is whether or not $B_d$ is divisible for
$d= 5.$ It is very likely that the answer is affirmative. However,
our proof does not seem to work in this case. The crucial
step in the proof of Theorems \ref{t2} is to check that the conditions of Lemma \ref{l1} is satisfied on the graph generated by the isometries.
Our proof in even dimension $d=2s$ is based on the fact
that if $O\in SO(2s)$ is a `generic' rotation then $O$ has no fixed point other than the
origin. Thus $T_b O$
has a fixed point for {\it every} vector $b\in \rr ^{2s},$
since $I-O$ is invertible, and $(I-O)\am (b)$ is a fixed point of $T_b O$.
This statement does not hold for dimension $d=2s+1$. Furthermore it can be easily shown that
1 is an eigenvalue of a `generic' rotation $O\in SO(2s+1)$ with multiplicity at least $1$.

However, it is also not clear if the method applied in \cite{KL} works for $d=5$. The result of \cite{KL} is based on the fact that if $O_0 \stb O_N \in SO(3)$ are `generic' rotations, $b\in \rr^3$
is a `generic' vector and $F=T_b O_0 ,$ then a nonempty reduced
word on the alphabet $O_1^{\pm 1} \stb O_N^{\pm 1}$ and $F^{\pm 1}$ has a fixed point only if the word is a conjugate of a word
on the alphabet $O_1^{\pm 1} \stb O_N^{\pm 1}.$ (See \cite[Lemma 3.5]{KL} for the precise statement.) Unfortunately, this
statement does not generalize for higher dimensions. The generalization has a difficulty.
The authors use the fact that the axis of a `generic' rotation $O$ can be expressed by the entries of the matrix $O$.
On the other hand, Borel \cite{B} proved that for every odd $d\ge 3$ there is a dense subset of pairs $(A, B)$ in $(SO(d))^2$ such that each pair generates a locally commutative group. Thus if two `generic' rotations have a common axis, then they commute. (This property is called the locally commutativity.)
Still, we conjecture that the corresponding graph has property \eqref{coin} in every dimension $d\ge3$.

\noi
{\tt Department of Stochastics, Faculty of Natural Sciences

\noi
Budapest University of Technology and Economics

\noi
and

\noi
MTA-BME Stochastics Research Group (04118)

\noi
Budapest, M\"uegyetem rkp. 3.

\noi
1111 Hungary

\msk \noi
{\rm e-mail:} kigergo57@gmail.com

\bsk
\noi
Department of Algebra

\noi
E\"otv\"os Lor\'and University

\noi
Budapest, P\'azm\'any P\'eter s\'et\'any 1/C

\noi
1117 Hungary

\msk \noi
{\rm e-mail:} gsomlai@cs.elte.hu

\end{document}